\def \version {2024--08--07}

\documentclass[12pt]{article}

\usepackage{amsmath}
\usepackage{amssymb}
\usepackage{amsthm}
\usepackage{MnSymbol,enumitem}
\usepackage{lineno}

\newtheorem{thm}{Theorem}[section]
\def \btm {\begin{thm}}
\def \etm {\end{thm}}
\newtheorem{prp}[thm]{Proposition}
\def \bpn {\begin{prp}}
\def \epn {\end{prp}}
\newtheorem{lem}[thm]{Lemma}
\def \blm {\begin{lem}}
\def \elm {\end{lem}}
\newtheorem{obs}[thm]{Observation}
\def \bob {\begin{obs}}
\def \eob {\end{obs}}
\newtheorem{rmk}[thm]{Remark}
\def \brm {\begin{rmk}}
\def \erm {\end{rmk}}
\newtheorem{cor}[thm]{Corollary}
\def \bcr {\begin{cor}}
\def \ecr {\end{cor}}
\newtheorem{con}[thm]{Conjecture}
\def \bcj {\begin{con}}
\def \ecj {\end{con}}
\newtheorem{prm}[thm]{Problem}
\def \bpm {\begin{prm}}
\def \epm {\end{prm}}
\newtheorem{dfn}[thm]{Definition}
\def \bdf {\begin{dfn}}
\def \edf {\end{dfn}}
\newtheorem{exa}[thm]{Example}
\def \bex {\begin{exa}}
\def \eex {\end{exa}}
\newtheorem{exas}[thm]{Examples}
\def \bexs {\begin{exas}}
\def \eexs {\end{exas}}

\def \bpf {\begin{proof}}
\def \epf {\end{proof}}
\def \appli {\paragraph{Application.}}

\def \bsk {\bigskip}
\def \msk {\medskip}
\def \ssk {\smallskip}
\def \nin {\noindent}
\def \smin {\setminus}
\def \es {\varnothing}

\def \cF {\mathcal{F}}

\def \arr {anti-Ramsey}
\def \Arr {Anti-Ramsey}
\def \socc {strong-odd-colored}
\def \cfcc {conflict-free-colored}
\def \cfcd {Cf-colored}
\def \lpcc {local-parity-colored}
\def \cpcc {class-parity-colored}
\def \spcc {strong-parity-colored}

\def \chiF {\chi_{_\cF}}
\def \AR {\mathrm{Ar}}
\def \lar {\mathrm{Lr}}
\def \oar {\mathrm{Od}}
\def \soar {\mathrm{Sod}}
\def \cfar {\mathrm{Cf}}
\def \lpar {\mathrm{Lp}}
\def \cpar {\mathrm{Cp}}
\def \parr {\mathrm{Sp}}
\newcommand{\Oar}[1]{\mathrm{\oar}(n,#1)}

\newcommand{\Ar}[1]{\AR(n,#1)}
\newcommand{\Lar}[1]{\lar(n,#1)}
\def \lkn {\lar(n,K_4)}

\def \lard {\Lar{K_4-e}}
\newcommand{\lod}[1]{\lar(#1,K_4-e)}
\newcommand{\Soar}[1]{\soar(n,#1)}
\newcommand{\Cfar}[1]{\mathrm{\cfar}(n,#1)}
\def \oard {\Soar{K_4-e}}
\newcommand{\sod}[1]{\soar(#1,K_4-e)}
\def \oapf {\Oar{P_5}}
\def \soapf {\Soar{P_5}}
\def \cfapf {\Cfar{P_5}}
\def \cfapf {\Cfar{P_5}}

\def \khp {K_3^+}

\def \ex {\mathrm{ex}}

\def \Lex {{\sf LEX}}
\def \Lexx {\mathrm{Lex}}
\def \UMIC {{\sf UMIC}}
\newcommand{\klek}[1]{{\sf #1-LEX}} 
\newcommand{\klex}[2]{{\sf (#1,#2)-LEX}} 
\newcommand{\krs}[1]{{\sf #1-RS}} 
\def \splc {{\sf SPLC}} 
\newcommand{\kcc}[1]{{\sf #1-CC}} 
\newcommand{\krum}[1]{{\sf #1-RUM}} 

\def \np {$\mathcal{NP}$}
\def \womo {odd-even ordering}

\def \erd {Erd\H os}
\def \sos {S\'os}
\def \tur {Tur\'an}

\begin{document}

\title{Monochromatic graph decompositions\\
 inspired by anti-Ramsey theory\\
  and the odd-coloring problem}
\author{Yair Caro\,\thanks{~Department of Mathematics, University of Haifa-Oranim, Tivon 36006, Israel}
 \and Zsolt Tuza\,\thanks{~HUN-REN Alfr\'ed R\'enyi Institute of Mathematics, Budapest, Hungary} $^,$\thanks{~Department of Computer Science and Systems Technology, University
of Pannonia, Veszpr\'em, Hungary}}
\date{\small Latest update on \version}
\maketitle

\begin{abstract}
We consider extremal edge-coloring problems inspired by the theory of
 \arr\ rainbow coloring, and further by odd-colorings and conflict-free colorings.

Let $G$ be a graph, and $\cF$ any given family of graphs.
For every integer $n \geq |G|$,
 let $f(n,G|\cF)$ denote the smallest integer $k$ such that any edge coloring of $K_n$ with at least $k$ colors forces a copy of $G$ in which each color class induces a member of $\cF$.  
Observe that in anti-Ramsey problems each color class is a single edge; i.e., $\cF=\{K_2\}$.

In our previous paper [arXiv:2405.19812], attention was given most\-ly to the
 case where $\cF$ is hereditary under subgraph inclusion. 
In the present work we consider coloring problems inspired by
 odd-coloring and conflict-free coloring.
As we shall see, dealing with these problems requires distinct
 additional tools to those used in our first paper on the
  subject.

Among the many results introduced in this paper, we mention:

\ssk

(1)\quad For every graph $G$, there exists a constant $c=c(G)$
 such that in any edge coloring of $K_n$ with at least $cn$
 colors there is a copy of $G$ in which every vertex $v$ is
 incident with an edge whose color appears only once among all
 edges incident with $v$.

\ssk

(2)\quad In sharp contrast to the above result we prove that if
 $\cF$ is the class of all odd graphs (having vertices with
 odd degrees only) then
  $f(n,K_k|\cF)= (1+o(1))\,\ex(n,K_{\lceil k/2 \rceil})$,
 which is quadratic for $k \geq 5$.

\ssk

(3)\quad We exactly determine $f(n,G|\cF)$ for small graphs
 when $\cF$ belongs to several families representing various
 odd/even coloring constraints. 

\bsk

\nin

\qquad\qquad \ \
parity colorings.

\bsk

\nin
\end{abstract}

\tableofcontents

\newpage

\section{Introduction}
\label{s:intro}

In this extensive work we introduce a large number of new
 functions related to the \Arr\ / Rainbow Theory of graphs,
 and present a first detailed study of them.
Doing so we continue the approach that we presented in our
 preceding paper \cite{ar-G}.
As a generalization of rainbow coloring, each monochromatic
 class is allowed to form a graph that belongs to a prescribed
 family $\cF$ of graphs, rather than required to be just
 a single edge.

In \cite{ar-G} we mainly considered hereditary families of graphs.
Here we focus on families, mostly not hereditary ones, that
 have parity conditions imposed on the degrees of the graphs
 allowed to form the color classes.

\subsection{Brief history of \arr\ theory}
\label{ss:history}

The rainbow coloring or \arr\ theorems date back to the work of
 \erd, Simonovits and \sos\ \cite{ESS-75}.
The main problem is: Given a graph $G$, at least how many colors
 are needed so that every edge coloring of $K_n$ with that many
 colors forces a copy of $G$ with all edges of it
 getting distinct colors?
This minimum is denoted by $\AR(n,G)$.
The main result of \cite{ESS-75} is:
Let $\chi_e(G) = \min \{ \chi(G - e) : e \in E(G) \}$,
 where $\chi$ denotes the chromatic number.
If $\chi_e(G) = k \geq 3$, then
 $\AR(n,G) = (1+o(1))\,\ex(n,K_k )$, and
 in particular $\AR(n,K_{k+1}) = \ex(n,K_k)+2$ for large $n$, 
 where $\ex(n,G)$ is the famous \tur\ number.
The proof relies heavily on the \erd--Stone--Simonovits theorem
 \cite{ES-46,ES-66} that states $\ex(n,G) = (1+o(1))\,\ex(n,K_k)$
 whenever $\chi(G) = k \geq 3$.
Decades later the tight formula for $\AR(n,K_{k+1})$ was
 extended also to small $n \geq k \geq 3$ \cite{S-04},
 moreover $\AR(n,C_k)$ was determined
 up to an additive constant \cite{MBNL-05}.
The subject is still very active today; see the survey \cite{dyn-surv}
 and some recent papers
  \cite{AC-24,GH-17,GH-22+,H-24t,WBZL-24,XLL-21,Y-21+}.

Recently in \cite{ar-G} the authors initiated the study of the
 following more general problem.
Let $G$ be a graph, and $\cF$ any given family of graphs.
For every integer $n \geq |G|$, let $f(n,G|\cF)$ denote the
 smallest integer $k$ such that any edge coloring of $K_n$
 with at least $k$ colors forces a copy of $G$ in which each
 color class induces a member of $\cF$.
Observe that in \arr\ problems each color class is a single
 edge; i.e., $\cF = \{ K_2 \}$.
A major result in \cite{ar-G} is:
Let $\cF$ be a hereditary family and let
 $\chiF(G) = \min \{ \chi(G - D)  :
  D \subset G ,\, D \in \cF \}$.
If  $\chiF(G) = k \geq 3$, then
 $f(n,G|\cF ) = (1+o(1))\,\ex(n,K_k )$.
Again a heavy use of the \erd--Stone--Simonovits theorem is
 inevitable, together with new tools of interest in their
 own right, e.g.\ the ``Independent Transversal Lemma'' in
  directed graphs of bounded outdegree.     
Also it is proved that if $G$ is stable with respect to
 $\cF$, namely $\chiF(G) = \chi(G) = k \geq 3$,
 then (regardless of whether or not $\cF$ is hereditary),
 $f(n,G|F) = (1+o(1))\,\ex(n,K_k )$ is valid.
Many examples of interesting and natural hereditary families and
 the implied results concerning $f(n,G|\cF)$ or $f(n,K_p |\cF )$
 are given in \cite{ar-G}.

Already in \cite{ar-G} we announced that in parallel we consider
 other similar problems inspired by odd-coloring and
 conflict-free coloring; cf.\ \cite{CPS-22} and \cite{CPS-23},
  respectively, with the several follow-ups and the references therein.
This track of research emerged from a
 conversation with Riste \v Skrekovski, and also is
 motivated by the famous theorem of Pyber \cite{P-91}
 stating that every graph has an edge decomposition into
 at most four odd subgraphs (and every multigraph without loops
 into at most six).
\Arr-type problems related to conflict-free nd odd-colorings
 constitute the main subject of the current paper.

\subsection{Hierarchy of some basic invariants under parity constraints}
\label{ss:hierarchy}

Beside the classical \arr\ numbers $\AR(n,G)$ we define seven
 further notions.
Their hierarchic relations are exhibited in Table \ref{tab:hierarc}.

\bdf
Let\/ $\psi$ be an edge coloring of\/ $K_n$, and let\/ $G$ be
 a given graph.
A subgraph\/ $H\cong G$ of\/ $K_n^\psi$ is called
 \begin{itemize}
  \item \emph{rainbow} if all its edges have mutually distinct colors;
  \item \emph{proper} or\/ \emph{local rainbow} if\/ $\psi$ induces a proper edge coloring on\/ $H$;
  \item \emph{\socc}, or just\/ \emph{strong}, if
    each color class induces an odd graph;
  \item \emph{odd-colored}, or shortly\/ \emph{weak} (as opposed to
   ``strong'') if at each vertex at least one color
   occurs on an odd number of edges in\/ $H$;
  \item \emph{\cfcc}---\/\emph{\cfcd}, for short---if at each vertex
   at least one color occurs on exactly one edge of\/ $H$;
  \item \emph{\spcc} if all color classes induce odd graphs
    or all color classes induce even graphs;
  \item \emph{\cpcc} if the edges of each color\/ $c$ form
    either an odd graph  or an even graph (but distinct
    color classes are not required to have the same parity);
  \item \emph{\lpcc} if, at each vertex, every incident color class
   has an odd number of edges or every incident color class
   has an even number of edges (but for distinct vertices
   the parity may not be the same).
 \end{itemize}
We denote by\/ $\AR(n,G)$ /\/ $\lar(n,G)$ /\/ $\soar(n,G)$ /\/
 $\oar(n,G)$ /\/ $\cfar(n,G)$ /\/
  $\parr(n,G)$ /\/ $\cpar(n,G)$ /\/ $\lpar(n,G)$
 the smallest $m$ such that, for every\/ $\psi$ with at least\/
 $m$ colors,\/ $K_n^\psi$ contains a rainbow / proper / \socc\ /
 odd-colored / \cfcd\ / \spcc\ / \cpcc\ / \lpcc\ subgraph
 isomorphic to\/ $G$, respectively.
\edf

For easier comparison among the definitions, we include
 two further tables from different aspects.
Table \ref{tab:fn} summarizes the
 combinations of local conditions assumed for each vertex,
 from which $\phi(n,G)$ is derived for the four functions
 $\phi\in\{\oar,\soar,\cfar,\lar\}$.
Table \ref{tab:loc-glob} exhibits the differences between
 assumptions posed globally for all vertices or
 locally at each vertex, involving the five functions
 $\phi\in\{\AR,\lar,\parr,\lpar,\cpar\}$.

\renewcommand{\arraystretch}{1.5}
\begin{center}
\begin{table}[ht]
\begin{center}
{\small
\begin{tabular}{ccccccc}
\hline
 &&  && \textbf{AR} && \\
 &&  && $\Downarrow$ && \\
 &&  && \textbf{LR} && \\
 &&  & $\Swarrow$ && $\Searrow$ & \\
 &&   \textbf{SOD} &&&& CF \\
 & $\Swarrow$ && $\Searrow$ && $\Swarrow$ & \\
 \textbf{SP} &&  && OD && \\
  $\Downarrow$ & $\Searrow$ &&&&& \\
 \textbf{CP} && \textbf{LP} &&&& \\
\hline
\end{tabular}
}
\end{center}
\caption{Hierarchy of eight subclasses of \arr\ colorings;
 boldface indicates possible quadratic growth, the other two
 classes are linearly bounded, as we shall see in the sequel.
   \label{tab:hierarc}}
\end{table}
\end{center}

\renewcommand{\arraystretch}{1.5}
\begin{center}
\begin{table}[ht]
\begin{center}
\begin{tabular}{cc|cccc}
color degree &&& odd && 1 \\
\hline
\hline
some color &&& $\oar$ && $\cfar$ \\
all colors &&& $\soar$ && $\lar$
\end{tabular}
\end{center}
\caption{Conditions for the four local variants of \arr\ functions. \label{tab:fn} }
\end{table}
\end{center}

\renewcommand{\arraystretch}{1.5}
\begin{center}
\begin{table}[ht]
\begin{center}
\begin{tabular}{cc|cccccc}
condition &&& global && local && color class \\
\hline
\hline
rainbow &&& $\AR$ && $\lar$ && --- \\
same parity &&& $\parr$ && $\lpar$ && $\cpar$
\end{tabular}
\end{center}
\caption{Global vs.\ local conditions on \arr\ functions. \label{tab:loc-glob} }
\end{table}
\end{center}

\paragraph{Representation in terms of $f(n,G|\cF)$.}

Four of the eight functions above,
 namely $\cfar(n,G)$, $\oar(n,G)$, $\parr(n,G)$, and $\lpar(n,G)$,
 are not possible to write in an equivalent form of
 $f(n,G|\cF)$ with any family $\cF$ of graphs;
two of them can be characterized with hereditary families:
 \begin{center}
  $\AR(n,G)$ $\longrightarrow$ $\cF = \{K_2\}$ (single edge), \\
  $\lar(n,G)$ $\longrightarrow$ $\cF = \{tK_2 : t\geq 1\}$ (matchings); \\
 \end{center}
and two of them with non-hereditary families:
 \begin{center}
  $\soar(n,G)$ $\longrightarrow$ family of all odd graphs, \\
  $\cpar(n,G)$ $\longrightarrow$ all odd graphs and all even graphs.
 \end{center}

\subsection{Summary of results}
\label{ss:summary}

The structure of this paper is explained in the Contents section;
 here we give just a brief list of few representative results.
We first emphasize the heavy use of \erd--Rado
 canonical coloring theorem (Theorem \ref{t:ERCT})
 along many results in this paper.
It seems that this theorem was rarely, if at all, used in
 anti-Ramsey theory previously, while in the present paper,
 introducing the parity-dependent anti-Ramsey parameters,
 this theorem has become rather useful. 

Clearly it is impossible to summarize in a short list
 all the many results obtained, hence we shall bring a couple
 of them which we hope represent the flavor of this work. 

From Section \ref{s:general} we choose to state the following
 theorem, showing that six of the eight parameters have
 quadratic growth in $n$, as indicated in the hierarchy
 presented in Table \ref{tab:hierarc}. 
For comparison we recall the tight formula
 $\AR(n,K_p) = \ex(n,K_{p-1}) +2$ for $p \geq 4$,
 proved in \cite{S-04}, that has a leading coefficient
 substantially larger than the one below.   

Theorem \ref{t:quad-asymp}: For every $p \geq 5$, all of
 $\lar(n,K_p)$, $\soar(n,K_p)$, $\parr(n,K_p)$,
 $\cpar(n, K_p)$, $\lpar(n, K_p)$ are asymptotically equal to
 $(1+o(1))\,\ex(n, K_{\lceil p/2 \rceil})$.  

From section 3 we choose a theorem that dictates the values for stars whenever $n$ is not too small. 

Theorem \ref{t:star}:
For $n \geq 2r$ we have $\soar(n,K_{1,r}) = 1$ if
 $r \equiv 1$ {\rm (mod 2)}, and $\soar(n,K_{1,r}) = 2$
 if\/ $r \equiv 0$ {\rm (mod 2)}.
Moreover, the condition\/ $n\geq 2r$ is tight in both cases.

Another result from Section \ref{s:degrees} shows a dichotomy
 in the behavior of three parameters:

Theorem \ref{t:odd}: If all vertices have odd degree in $G$,
 then $\soar(n,G)$, $\parr(n, G)$, $\lpar(n,G)$
 either all tend to infinity with $n$,
 or all are equal to 1 for every $n$ sufficiently large. 

A major theorem of Section \ref{s:CF-linear}, completing the support of Table \ref{tab:hierarc} in the branch of linear growth, is: 

Theorem \ref{t:CF-G_p}: If $G$ is a graph on $p$ non-isolated
 vertices, then
  $\oar(n,G) \leq \cfar(n,G) \leq (p-2)n -
   \lfloor p^2\!/2 \rfloor +p+1$. 

In Section \ref{s:more-small} we mostly deal with
 small graphs up to 4 edges as well as the graphs
 $P_6$, $K_4-e$, and $K_4$.
These results are summarized in Table \ref{tab:rsmall}
 on page \pageref{tab:rsmall}.  
For $P_6$ and all graphs with at most four edges we have exactly
 determined the values of all the newly introduced parameters.
However, several cases are left open regarding $K_4-e$ and $K_4$. 
Yet we mention here one result concerning $K_4$ where
 the order of growth is determined. 

Theorem \ref{t:l-K4}:
 We have $\lar(n, K_4) = \Theta(n^{3/2})$ as $n \to\infty$. 

Lastly, it is inevitable that from such a collection of
 parameters lots of new problems will emerge.
Section \ref{s:conclude} collects many of them for future research.

\subsection{Standard definitions and notation}

Throughout the paper, we consider simple undirected graphs
 $G$, without loops and multiple edges,
 with vertex set $V(G)$ and edge set $E(G)$.
We write $|G|$ for the order $|V(G)|$ of $G$.
The degree of vertex $v$ in graph $G$ is denoted by $d_G(v)$,
   abbreviated as $d(v)$ when $G$ is understood.
Also, as usual, $\delta(G)$ and $\Delta(G)$ denote the minimum
 and maximum degree in $G$, respectively.

We say that $G$ is an \emph{odd graph} if the degree of all its
 vertices is odd; and $G$ is an \emph{even graph} if all
 vertex degrees in $G$ are even.

Particular types of graphs are the path $P_n$, the cycle $C_n$,
 and the complete graph $K_n$, each of order $n$, and the
 complete bipartite graph $K_{p,q}$ with $p$ and $q$ vertices
 in its partite sets.
Where a considered edge coloring $\psi$ of $K_n$ has to be
 emphasized, we write $K_n^\psi$.

Some important graph operations are edge insertion $G+e$
 (with $e\subset V(G)$ and $e\notin E(G)$), edge deletion
 $G-e$ (with $e\in E(G)$), and the vertex-disjoint union
 $G_1\cup G_2$ of two graphs $G_1$ and $G_2$.
The vertex-disjoint union of $t$ copies of $G$ is denoted as $tG$.

A famous extremal graph-theoretic function of great importance
 in the current context is the \emph{\tur\ number} $\ex(n,F)$
 of a ``forbidden graph'' $F$.
It is defined as the maximum number of edges in a graph of
 order $n$ that contains no subgraph isomorphic to $F$.

\section{Basic coloring patterns, general inequalities, and graph operations}
\label{s:general}

In this section we present three kinds of material.
The first part describes explicit constructions of coloring
 patterns that can be used for lower bounds on \arr\ numbers.
The second part deals with inequalities based on hierarchy
 among notions and simple structural observations.
The third part disccusses the effect of some operations
 that simplify the  determination of \arr\ parameters,
 pointing out reducibility relations among them.

\subsection{Coloring patterns}

Let us categorize the edge coloring patters of $K_n$ below in
 two major types: homogeneous and compound.

Homogeneous ones are fundamental and will be useful in proving
 upper bounds on several \arr\ functions under study.
On the other hand,
 the various constructions of compound patterns are tools for
 proving lower bounds.
Already after the definitions, we will indicate several ways
 they can be used in that direction.

\subsection{Homogeneous coloring patterns}

The three fundamental types of homogeneous coloring patterns are:
 \begin{itemize}
  \item monochromatic---all edges of $K_n$ are assigned with
    the same color;
  \item rainbow---each edge of $K_n$ has its private color,
    distinct from all the colors of the other edges;
  \item \Lex\ coloring---labeling the vertices of $K_n$ as
    $v_1,v_2,\dots,v_n$, for all $1\leq i<j\leq n$ the edge
    $v_iv_j$ is assigned with color $j-1$.
 \end{itemize}
Hence the classical \Lex\ coloring\footnote{In
   some papers \Lex\ is also termed \UMIC\
 as an abbreviation for Universal Majority Index Coloring;
 but throughout this work we use the more standard name \Lex.}
 of $K_n$ assumes a sequential order on the vertex set, and
 partitions the edge set into $n-1$ color classes, which are stars.

The following important Ramsey-type result connects the above
 three coloring patterns.

\btm {\bf\itshape (\erd--Rado Canonical Theorem \cite{ER-50})}   \label{t:ERCT}
For every\/ $k\geq 3$ there exists an integer\/ $n_0(k)$ such
 that, for every\/ $n\geq n_0(k)$, in every edge coloring of\/
 $K_n$, some\/ $k$ vertices induce a copy of\/ $K_k$ which is
 either monochromatic or rainbow or\/ \Lex-colored.
\etm

\subsection{Compound coloring patterns}

The three basic patterns introduced above can be combined with
 each other in many ways; here we describe some of those possibilities.
The ones used in the sequel are collected in Table \ref{tab:colpatt},
 the others partly provide tools for lower bounds in
 \cite{ar-G} and partly can be applied to improve estimates
 that are asymptotically tight but better error terms can be achieved.

\begin{center}
 \begin{table}
\begin{center}
  \begin{tabular}{ccc}
  \hline
   pattern combination && applied first in \\
  \hline
   rainbow $K_p$, monochromatic $K_n-K_p$ && Proposition \ref{p:num-even} \\
   \klek{h}, rainbow $K_h$ followed by \Lex && Theorem \ref{t:even} \\
   \krs{h}, rainbow $K_n-K_{n-h}$, monochromatic $K_{n-h}$ && Theorem \ref{t:path-cyc-LB} \\
   \Lex\ with rainbow perfect matching && Theorem \ref{t:l-K4-e} \\
   rainbow spanning graph, monochromatic complement && Theorem \ref{t:l-K4} \\
   \Lex\ with rainbow spanning star && Theorem \ref{t:s-K4} \\
  \hline
  \end{tabular}
\end{center}
   \caption{Compound coloring patterns where the number of
    colors grows as a function of $n$.}   \label{tab:colpatt}
 \end{table}
\end{center}

\subsubsection{Generalizations of \Lex\ coloring}
\label{ss:lex}

Recall that \Lex\ assumes a sequential order $v_1,v_2,\dots,v_n$
 on the vertices of $K_n$.
We generalize \Lex\ to a pattern with two parameters.
\label{mod:lex}

Let $k \geq 1$ and $h \geq k+1$.
Compose the coloring \klex{k}{h} as follows.
 \begin{itemize}
  \item For each $j=h+1,\dots,n$ and all $k\leq l<j$ assign
   the color $\psi(v_lv_j)=j-h$.
   Those colors range from 1 to $n-h$.
  \item For each $j=h+1,\dots,n$ and each $1\leq l<k$ assign
   a private color $\psi(v_lv_j)$, ranging from $n-h+1$
   to $(n-h)k$.
   (This step is void if $k=1$.)
  \item Take a rainbow $K_h$ on the vertices $v_1,\dots,v_h$
   using the colors $(n-h)k + 1,\dots,(n-h)k +  h(h-1)/2$.
 \end{itemize}
Here \Lex\ is obtained by putting $k = 1$ and  $h=1$ or $h = 2$.
Fixing $k=1$, the intermediate one-parameter case with $h>2$,
 termed \klek{h},
 is also of interest; the number of colors in that pattern is
  $$
    \Lexx(n,h):=\Lexx(n,1,h) = n-h + h(h-1)/2 = n + h(h-3)/2.
  $$
In general, the number of colors is
 $$
   \Lexx(n,k,h) =  (n-h)k +  h(h-1)/2 \,.
 $$

\appli

If $|G|=h+1$ and $\delta(G) =  k+1$, then
 $$
   \lar(n,G)  \geq \Lexx(n,k,h) +1  = (n-h)k +  h(h-1)/2 +1.
 $$
Indeed, taking $K_n^\psi$ with the \klex{k}{h}\ coloring,
 in any copy of $G$ the vertex $v_j$ with highest index
 satisfies $j>h$, therefore only $k$ colors occur on the
 edges from $v_j$ to its neighbors (all with smaller indices),
 but $d(v_j) \geq k+1$, hence some color appears at least
 twice at $v_j$ and so the coloring of $G$ cannot be proper.
This coloring sometimes provides a useful lower bound when
 $\chi(G) \leq 4$, as $\lar(n,G)$ is determined by the
 hereditary family $\cF = \{ tK_2 : t \geq 1\}$, and in view of
 the results presented in the introduction from \cite{ar-G}.

\subsubsection{Split coloring combined with Lex}

We present here two variants.
Their common feature is that the vertex set of $K_n$ is
 split into two parts, $V(K_n)=A \cup B$, where $|A| = k$
 and $|B| = n - k$ for a given $k\geq 1$, and
 \begin{itemize}
  \item on the edges of the clique on $A$ and all edges
   from $A$ to $B$, all colors are distinct.
 \end{itemize}

\nin
$(i)$\quad
In the coloring \krs{k}\ ($k$-rainbow split graph coloring)
 \begin{itemize}
  \item one extra color is used on all the other edges,
   to make the clique on $B$ monochromatic.
 \end{itemize}
The number of colors used is
 $$
   s(n,k) =  k(n-k) + \binom{k}{2} +1 = kn - \binom{k+1}{2} +1.
 $$

\appli

If $\delta(G) \geq k +2$, then
 $$
   \lar(n,G) \geq s(n,k) +1  = kn - \binom{k+1}{2} +2
 $$
 because any copy of $G$ has a vertex $v$ in $B$,
 which is incident with no more than $k$ distinct colors
  towards $A$ and just one color in $B$,
 hence at least two edges at $v$ must have the same color
 inside $B$.
Consequently, $G$ is not properly colored under \krs{k}.
This coloring is sometimes useful for lower bounds when $\chi(G)\leq 4$.

\bsk

\nin
$(ii)$\quad
In the Split-Lex coloring pattern \splc\ we again take the
 rainbow set of all edges meeting $A$, and
 \begin{itemize}
  \item apply the \klex{t}{h} coloring inside $B$.
 \end{itemize}
The number of colors used is
\begin{eqnarray}
   s(n, k, t, h) &=& s(n,k) -1  +
    \Lexx(n-k,t,h) \nonumber \\
    &=&  k(n-k) + \binom{k}{2} + (n-k-h)t + \binom{h}{2}. \nonumber 
\end{eqnarray}
Already the very particular case $k=t=1$, $h=2$
 (rainbow spanning star and \Lex\ on its leaves) with
 $s(n, k, t, h)=2n-3$ is of interest.

\appli  

This pattern provides the currently best known lower bound
 on the strong odd \arr\ number of $K_4$, that is
 $\soar(n,K_4) \geq 2n -2$.

\subsubsection{Clique coloring}

Let $k\geq 2$ be given.
To obtain the $k$-Clique coloring pattern \kcc{k}, we write
 $n$ in the form $n = qk + r$, and partition the vertex set
 of $K_n$ as $V(K_n)=A_0\cup A_1\cup \cdots \cup A_q$,
 where $|A_0|=r<k$ (possibly $A_0=\es$)
  and $|A_1| = \ldots = |A_q|=k$.  
 \begin{itemize}
  \item Inside every $A_i$ ($i=0,1,\dots,q$) let the complete
   subgraph get the rainbow coloring, each color appearing
   in just one $A_i$.
  \item Assign a fresh new color to all other edges, to make
   a monochromatic complete multipartite graph.
 \end{itemize}
The number of color used is
 $$
   q(n,k) = (n-r)(k-1)/2 + \binom{r}{2} +1.
 $$

\appli

This coloring is sometimes useful in giving lower bounds when
 $\chi(G) \leq 4$, and also can be applied to derive a
 lower bound on $\lar(n,G)$ in terms of maximum degree;
  cf.\ Proposition \ref{p:Deg-LB}.

\subsubsection{Rainbow multipartite coloring}
\label{ss:krum}

Let $k\geq 2$ be given.
We present two versions of a pattern, a simpler and a
 more involved one.
For both of them we take a balanced $k$-partition of the
 vertex set, $V(K_n)=A_1\cup \cdots \cup A_k$ with
  $\lfloor n/k \rfloor \leq |A_i| \leq \lceil n/k \rceil$
 for all $1\leq i\leq k$.
\begin{itemize}
 \item Assign mutually distinct colors to all edges
  between distinct parts $A_i,A_j$.
\end{itemize}

\nin
$(i)$\quad
In the simpler form, \krum{k} coloring,
\begin{itemize}
 \item apply \Lex\ inside each $A_i$, no color appearing
  in more than one part.
\end{itemize}
The balanced multipartite graph uses $\ex(n,K_{k+1})$ colors,
 and inside each part $A_i$ the number of colors is $|A_i|-1$.
Hence the total number of colors is
 $$
   r(n,k) =  n-k + \ex(n,K_{k+1}).
 $$

\appli

This construction yields a quadratic lower bound on all
 boldface functions as indicated in Table \ref{tab:hierarc}.
More explicitly, the following asymptotics can be proved for
 five of the six functions; the exception is $\AR(n,K_p)$ whose
 behavior is substantially different as proved in \cite{ESS-75}.

\btm   \label{t:quad-asymp}
For every\/ $p \geq 5$, all of\/ $\lpar(n,K_p)$, $\cpar(n,K_p)$,
 $\parr(n,K_p)$, $\soar(n,K_p)$, $\lar(n,K_p)$ are asymptotically
 equal to\/ $(1+o(1))\,\ex(n,K_{\lceil p/2 \rceil})$.
\etm

\bpf
Consider the \krum{($\lceil$p/2$\rceil$--1)} coloring on $K_n$,
 which uses
  $\ex(n,K_{\lceil p/2 \rceil}) + n - \lceil p/2 \rceil + 1$
  colors.
Then any subgraph $K\cong K_p$ of $K_n$ contains at least three
 vertices in some $A_i$.
In the \Lex\ ordering the edges from the third vertex to the
 first and second vertices form a color class $P_3$; moreover,
 from the third vertex there is an edge either to another $A_j$
 or a vertex of higher index inside $A_i$.
The monochromatic $P_3$ is not allowed in a class parity coloring,
 and the presence of degrees 1 and 2 simultaneously at a vertex
 is not allowed in a local parity coloring.
This yields the lower bound
 $\ex(n,K_{\lceil p/2 \rceil})-\lceil p/2 \rceil+1$
 for all the five functions under consideration, because
 $\cpar$ and $\lpar$ are lower bounds on all of $\parr$,
 $\soar$, and $\lar$.

The matching asymptotic upper bound for $\lar(n,K_p)$ holds by
 Theorem 5.2 (i) of \cite{ar-G}, and it dominates the other four
 functions involved in the assertion.
\epf

Beyond the present setting, this pattern is among the few standard
 colorings that are heavily used in our paper \cite{ar-nonham}
 where we consider the problem of determining $f(n,G|\cF)$
 in full generality.
The families $\cF$ there neither are
 required to be hereditary, nor have they any relation to
 odd-coloring or conflict-free coloring.

\bsk

\nin
$(ii)$\quad
In the more complex form, rainbow coloring is combined with
 the generalization of \Lex:
\begin{itemize}
 \item apply \klex{t}{h}\ inside each $A_i$, no color
  appearing in more than one part.
\end{itemize}
Disregarding small deviation due to integer parts,
 the number of colors is approximately
 $$
   r(n,k,t,h) = \ex(n,K_{k+1}) + k \Lexx(n/k,t,h) =
     \ex(n,K_{k+1}) + k((n/k -h)t +  h(h-1)/2 )
 $$
that means a somewhat larger linear addition in the number of
 colors to $\ex(n,K_{k+1})$.

\appli

It can be show that this pattern gives better lower bounds for
 $\lar(n,K_p)$ than the simpler version.
Although for $p \geq 5$  the asymptotics
  $\lar(n,K_p) = (1+o(1))\,\ex(n,K_{\lceil p/2 \rceil})$ is
  well determined, this improvement may be relevant when
  exact results are concerned.

\subsection{Inequalities between \arr\ functions}

Some basic relations among the five functions not involving ``modulo 2'' constraints are collected next.

\bob   \label{ob:ineq}
For any graph\/ $G$ and any\/ $n\geq |G|$ we have:
 \begin{enumerate}
  \item 
 $\oar(n,G) \leq \soar(n,G) \leq	 \lar(n,G) \leq \AR(n,G)$;
  \item 
 $\oar(n,G) \leq \cfar(n,G) \leq \lar(n,G) \leq \AR(n,G)$;
  \item 
 the equalities\/ $\oar(n,G) =  \soar(n,G) = \cfar(n,G)  = \lar(n,G)$  are valid
 whenever\/ $G$ has maximum degree at most\/ $2$.
 \end{enumerate}
\eob

Concerning \arr\ functions involving parity, the following
 facts are valid.

\bpn   \label{p:parity}
Let\/ $G$ be any graph, and\/ $n\geq |G|$.
\begin{enumerate}
 \item For every\/ $G$ we have\/ $\soar(n,G) \geq \parr(n,G)$,
  moreover\/ $\parr(n,G) \geq \lpar(n,G)$  and\/
   $\parr(n,G) \geq \cpar(n,G)$.
 \item If there is a vertex of odd degree in\/ $G$, then\/
  $\parr(n,G) = \soar(n,G)$; and if\/ $G$ is an odd graph,
  then\/ $\lpar(n,G) = \parr(n,G) = \soar(n,G)$.
 \item If\/ $G$ has a component\/ $H$ with all degrees even, then\/
  \hbox{$\parr(n,G) \geq \cpar(n,G) \geq n$}
   is forced by\/ \Lex\ for\/ $H=C_3$ and\/
   \hbox{$\parr(n,G) \geq \cpar(n,G) \geq n +1$}
    by\/ \klek{3}\ for any other\/ $H$.
 \item If\/ $\Delta(G) \leq 2$, then\/ $\lpar(n,G) = 1$.
 \item If\/ $G$ is a linear forest (i.e., contains no cycles
  and\/ $\Delta(G) \leq 2$), then\/
   $\lar(n,G) = \cfar(n,G) = \oar(n,G) = \soar(n,G) =
    \parr(n,G) = \cpar(n,G)$.
 \item The\/ $\lpar$ condition is satisfied by
  every monochromatic graph and also by every rainbow graph.
\end{enumerate}
\epn

\bpf
\ 1.\quad
Strong odd coloring is the restricted version of strong parity
 coloring where even degrees are not allowed in any color class.
On the other hand, as compared to strong parity coloring,
 the color classes in a class parity coloring may mix odd graphs
 and even graphs, whereas in a local parity coloring also graps
 are allowed to be color classes that contain vertices with
 degrees from both parities.

\msk

\nin 2.\quad
At an odd-degree vertex, at least one color occurs an odd
 number of times, in any edge coloring.
This forces all colors to have odd degrees at all vertices,
 in every strong parity coloring.
In a general graph, the parity may vary vertex by vertex, but
 in odd graphs all parities must be odd, hence the difference
 between local parity and strong parity disappears.

\msk

\nin 3.\quad
In \Lex, using $n-1$ colors, the three edges of every triangle
 have exactly two colors that do not occur anywhere else in $G$,
 that is not a class parity coloring.

The pattern \klek{3} uses $n$ colors.
If the component $H$ of $G$ is an even graph other than $C_3$,
 then necessarily $|H|>3$.
Thus, in any copy of $G$ the edges at the vertex of highest
 index in $H$ form a monochromatic star of even degree,
 which is not allowed in a class parity coloring.

\msk

\nin 4.\quad
The edges at a vertex of degree 2 either are monochromatic or
 have color distribution $1+1$.
Both cases are allowed in local parity coloring.

\msk

\nin 5.\quad
All the six parameters listed in the assertion require that the
 color class present at a vertex of degree 1 be a single edge.
Induction yields that every allowed coloring of $G$ is a
 proper edge coloring.

\msk

\nin 6.\quad
Either there is only one color class at a vertex---having
 degree parity $d(v)$~mod~2---or all colors present at $v$
 occur just once, thus all are odd.
\epf

\subsection{Effect of graph operations}

Based on the principle behind Part 3 of Observation \ref{ob:ineq},
 the following further inequalities can be derived.

\bpn   \label{p:deg2}
\textbf{(Adding Edge Lemma)}
\
\begin{itemize}
 \item[$(i)$]
Let\/ $v,w$ be nonadjacent vertices of even degrees in a
 graph\/ $G$, and let\/ $G^+$ be the graph obtained from\/
 $G$ by inserting the new edge\/ $vw$.
Then\/ $\oar(n, G^+) \leq \oar(n, G)$.
 \item[$(ii)$]
If\/ $G$ is a spanning subgraph of a graph\/ $H$, and for
 every vertex\/ $v$ either\/ $d_H(v)=d_G(v)$ or\/ $d_H(v)$
 is odd, then\/ $\oar(n,H)\leq \oar(n,G)$.
In particular, if\/ $\Delta(H)=2k+1$ and\/
 $\delta(H) \geq 2k -1$, and each vertex of degree\/ $2k-1$
 in\/ $G$ has degree\/ $2k-1$ or\/ $2k+1$ in\/ $H$, then\/
 $\oar(n,H)\leq \oar(n,G)$.
 \item[$(iii)$]
Let\/ $v,w$ be nonadjacent vertices of degree\/ $2$ in a graph\/ $G$, and
 let\/ $G^+$ be the graph obtained from\/ $G$ by inserting
 the new edge\/ $vw$.
Then\/ $\cfar(n,G^+)\leq \cfar(n,G)$.
\end{itemize}
\epn

\bpf
For $\phi=\oar$ in $(i)$--$(ii)$ and $\phi=\cfar$ in $(iii)$
 we consider any edge coloring $\psi$
 of $K_n$ with at least $\phi(n,G)$ colors.

\msk

\nin
$(i)$\quad
By assumption, a weak copy of $G$ exists under $\psi$.
No matter what the color of $\psi(vw)$ is, after the insertion
 of edge $vw$ the degrees of $v$ and $w$ become odd, and by
 parity an odd color must occur at each of them in $G^+$.

\msk

\nin
$(ii)$\quad
Also here a weak copy of $G$ exists under $\psi$.
Let $\psi$ be a coloring realizing $\oar(n,G)$,
 and consider $H$.
A vertex of even degree in $H$ remains with the
 original coloring and has at least two odd colors.
A vertex of odd degree in $H$
 must have an incident odd-degree color by parity.

\msk

\nin
$(iii)$\quad
In a Cf-colored graph $H\subset K_n^\psi$, $H\cong G$, both
 $v$ and $w$ in $H$ are incident with two distinct colors.
Inserting the edge $vw$, the color distributions at $v$ and $w$
 are modified from $1+1$ to $1+2$ or $1+1+1$, both satisfying
 the requirements for $H^+\cong G^+$ to be a Cf-colored graph.
\epf

\subsection{Large jumps when adding/deleting edges}
\label{ss:jump}

Here we list some examples illustrating the possible effects of
 edge insertion, causing large jumps in most of the parameters,
 or being non-monotone as in the second variable of $\oar(n,G)$
 even on graphs with $\Delta(G)\leq 2$
  (cf.\ Proposition \ref{p:deg2}).

\begin{enumerate}
 \item The sequence $P_4 \longrightarrow C_4 \longrightarrow
   K_4-e \longrightarrow K_4$ obtained by inserting edges
  between vertices of degree at most 2 has
   $\oar(n,P_4) = 3$,  $\oar(n,C_4) = n +1$,  $\oar(n,K_4-e) = 3$,  $\oar(n,K_4) = 1$.
  Hence any large increase and any large decrease can occur,
   even of linear order, despite that\/ $\oar(n,G)$ is linear
   for every graph\/ $G$.

  In comparison, the same sequence of graphs yields the
   following spectacular increase of the growth orders in $n$
   for $\AR(n,G)$, which clearly is monotone in terms of $G$\,:
\begin{center}
constant $\to$ linear $\to$ of order $n^{3/2}$ $\to$ quadratic.
\end{center}
 \item For cycles of even length\/ $k$ we know from
   Theorem \ref{t:path-cyc-LB} that\/ $\oar(n,C_k)$ grows with\/
    $\lfloor (k-1)/3 \rfloor n$ at least.
   However, inserting a perfect matching\/ $M$ we obtain an odd graph
    (more explicitly\/ $3$-regular) which has\/ $\oar=1$.
   It follows that inserting the edges of\/ $M$ one by one,
    if\/ $k\equiv 4 \ (\mathrm{mod} \ 6)$, the average fall
    of\/ $\oar$ per insertion can be estimated with nearly\/
     $2n/3$ from below.
   The worst case of decrease during the insertion of a single
    edge (between two vertices of degree\/ $2$, and also
    between any two vertices) is currently not known.
 \item From the results of \cite{ar-G} we know that
 $\lar(n,K_5-P_3)=o(n^2)$, caused by the fact that
  there is a way to omit $2K_2$ from $K_5-P_3$ to obtain
  the bipartite graph $K_{2,3}$.
 On the other hand, $\lar(n,K_5-e) = (1+o(1))\,\ex(n,K_3)$ because
 $K_3\subset (K_5-e) - 2K_2$, no matter how the removal of
 $2K_2$ is performed.
  Hence the insertion of a new edge into
 $K_5-P_3$ makes $\lar$ jump from subquadratic to quadratic.
 \item The parameters $\phi(n,C_5)$ are linear in $n$
 for all the eight functions
   $\phi\in\{\AR,\lar,\soar,\parr,\cpar,\lpar,\cfar,\oar\}$,
 and $\phi(n,K_5)$ is quadratic for six of them,
  except for $\cfar$ and $\oar$.
 Hence inserting edges one by one to reach $K_5$ from $C_5$,
  interesting jumps occur in those functions.
\end{enumerate}

\section{Vertex degrees}
\label{s:degrees}

\subsection{Some general inequalities}
\label{ss:deg-general}

\bpn   \label{p:Deg-LB}
Let\/ $G$ be a graph with maximum degree\/ $\Delta = \Delta(G)$.
Suppose\/ $n \equiv r \ (\mathrm{mod} \ \Delta -1)$. 
Then\/ $\lar(n,G) \geq q(n,\Delta-1) + 1
  =  ( n-r)(\Delta-2)/2 +\binom{r}{2} +2$.
\epn

\bpf
Consider $(\Delta-1)$-Clique coloring of $K_n$, with
 $q(n,\Delta-1)$ colors. 
A vertex $v$ of degree $\Delta$ in any copy of $G$
 can get only $\Delta-2$ distinct
 colors in a rainbow $K_{\Delta-1}$,
 hence at least two edges of the same extra color
 must be incident with $v$
  and therefore no properly colored $G$ occurs.
\epf

\bpn   \label{p:num-even}
If the number of vertices with positive even degree is\/
 $s>0$ in a graph\/ $G$, then\/ $\binom{s-1}{2}+2$ is a
 lower bound on all of\/ $\AR(n,G)$, $\lar(n,G)$,
 $\cfar(n,G)$, $\soar(n,G)$, and\/ $\oar(n,G)$.
\epn

\bpf
Based on Observation \ref{ob:ineq}, we only have to prove
 $\oar(n,G) \geq \binom{s-1}{2}+2$.
Inside $K_n$, take a rainbow $K_{s-1}$, and make $K_n-K_{s-1}$
 monochromatic in a new color.
At least one of the even-degree vertices belongs to
 $V(K_n)\smin V(K_{s-1})$, and its incident edges form a
 monochromatic star.
Hence, this coloring with $\binom{s-1}{2}+1$ colors does not
 admint any odd-colored copy of $G$.
\epf

It will be shown in Section \ref{ss:deg-odd-even} that a
 substantial improvement can be given if all vertices have
 even degrees.

\subsection{Stars}
\label{ss:star}

For stars we can present tight results under the assumption that
 the number $n$ of vertices is not very small.
In the next theorem we do not make efforts to optimize the
 bound on $n$.

\btm   \label{t:star}
\
 \begin{itemize}
  \item[$(i)$] If\/ $r \equiv 1$ {\rm (mod 2)} and\/
   $n \geq 2r$, then\/ $\soar(n,K_{1,r}) = 1$.
  \item[$(ii)$] If\/ $r \equiv 0$ {\rm (mod 2)} and\/
   $n \geq 2r$, then\/ $\soar(n,K_{1,r}) = 2$.
 \end{itemize}
Moreover, the condition\/ $n\geq 2r$ is tight in both cases.
\etm

\bpf
If $n=2r-1$, then we can decompose $K_n$ into $r-1$ Hamiltonian
 cycles as color classes.
Under this coloring a strong-odd-colored $K_{1,r}$ does not
 occur because only one edge can be selected from each color.
Hence $n\geq 2r$ is necessary to assume.

If $r$ is odd, then $K_{1,r}$ is an odd graph, and the
 monochromatic $K_n$ contains it for any $n\geq r+1$.
If $r$ is even, then $K_{1,r}$ itself is not allowed to be a
 color class in a strong odd coloring, hence $\soar(n,K_{1,r})>1$.
So, from now on we may restrict ourselves to colorings that
 use at least two colors.

If $n\geq 2r$, consider any non-monochromatic coloring of $K_n$.
Choose a vertex $v$ whose incident edges are also
 non-monochromatic.
If some color occurs at least $r$ times at $v$, then we find
 monochromatic $k_{1,r}$ if $r$ is odd, or monochromatic
 $k_{1,r-1}$ with a further edge from another color class
 if $r$ is even, and the proof is done.

Otherwise, select a star $S$ of any $2r-1$ edges at $v$.
Assume that color $i$ has $d_i$ edges in $S$, where the degrees
 $d_1\geq d_2\geq\ldots$ are in decreasing order (re-indexed
 if necessary).
Denote by $k$ the number of colors in $S$; we certainly have
 $k\geq 2$, because $d_1<r$.
We may assume $k<r$ (otherwise a rainbow $K_{1,r}$ is present
 and we have nothing to prove).
This assumption implies $d_1\geq 3$.

If $k \equiv r$ (mod 2), then we select one edge from each of the
 $k$ color classes, and sequentially supplement them with
 monochromatic pairs of edges, until at most one unselected edge
 remains in each class.
The process stops when $d_i$ or $d_i-1$ edges of color $i$ are
 selected, whichever is odd; in either case it means at least
 $d_i/2$ edges, and in fact at least $\frac{3}{4} d_i$
  unless $d_i=2$.
Thus, the number of selected edges is not smaller than
 $\lceil (d_1+\ldots+d_k)/2 \rceil = \lceil (2r-1)/2 \rceil = r$,
 and there is an intermediate step with exactly $r$ edges.
In that moment a strong-odd-colored $K_{1,r}$ is obtained.

If $k \not\equiv r$ (mod 2), we cannot use all colors;
 we then perform the above selection in the first $k-1$ colors.

The case $k=2$ forces $r$ to be odd.
Then a $K_{1,r}$ in color 1 is present, because
 $d_1\geq \lceil (d_1+d_2)/2 \rceil = \lceil (2r-1)/2 \rceil = r$.

If $k\geq 3$, recall that $d_1\geq 3$ and at least
 $\max\,(3,d_i-1)\geq\max\,(3,\lceil \frac{3}{4} d_i \rceil)$
 edges can be selected from each color class $i<k$
 of size $d_i\geq 3$ (and one edge if $d_i=2$).
If $d_k\leq 2$, we have $\max\,(3,d_i-1)\geq (d_1+d_k)/2$, so
 eventually at least half of the $2r-1$ edges will be selected.
And if $d_k\geq 3$, then $d_1+d_2\geq \frac{2}{3}(d_1+d_2+d_k)$,
 moreover all $d_i$ are at least 3.
Thus, we can apply
  $\frac{3}{4}(d_1+d_2)\geq \frac{1}{2}(d_1+d_2+d_k)$,
 completing the proof.
\epf

\bcr
$\oar(n, K_{1,r}) = \soar(n,K_{1,r})$ for\/ $n \geq 2r$.
\ecr

\bpf
Clearly, $1\leq \oar(n, K_{1,r}) \leq \soar(n,K_{1,r})$.
Hence any coloring satisfies the requirements
 for both parameters if $r$ is odd.
If $r$ is even, we also have $\oar(n, K_{1,r}) \geq 2$,
 that matches the upper bound.
\epf

\btm
If\/ $r \geq 3$ and\/ $n \geq 2r - 2$, then\/
 $\cfar(n, K_{1,r}) = 2$.
Moreover, the condition\/ $n \geq 2r  - 2$ is best possible.
\etm

\bpf
Clearly, we need at least two colors at the center of the star.
Consider any non-monochromatic coloring of $K_n$, and let $v$
 be a vertex incident with at least two colors.
The degree of $v$ is at least $2r-3$, hence omitting the smallest
 color class at $v$ there remain at least $r-1$ edges.
Thus, we can take $r-1$ of them together with one edge from the
 smallest color class, and obtain a Cf-colored $K_{1,r}$.

For $n = 2r-3$ the degree is $2r - 4$ and we can color the
 edges with two colors such that each color class is an
 $(r-2)$-regular spanning subgraph of $K_n$.
Then any $r$ edges at any vertex contain more than one edge
 from each color, hence no Cf-colored $K_{1,r}$ occurs.
\epf

\subsection{Shortest brooms}
\label{ss:broom}

\def \brmk {K_{1,k}^+}

We denote by $\brmk$ the tree obtained from the star $K_{1,k}$
 by attaching a pendant edge to one of its leaves.
It is the same as the double star $D_{1,k-1}$ whose two central
 vertices have degree 2 and $k$, respectively; it is a
 caterpillar, and also can be viewed as the broom graph obtained
 by identifying the central vertex of $K_{1,k-1}$ with an end of
 the path $P_3$.

Before turning to this graph, let us mention that the \arr\
 number of stars with $k\geq 3$ edges was determined in
  \cite{J-02,MB-06,C+al-23+} as
$$
  \AR(n,K_{1,k}) = \lfloor n(k - 2) / 2 \rfloor
   + \lfloor n / (n - k + 2) \rfloor + \varepsilon
$$
 where $\varepsilon\in\{0,1\}$.
For example, $\varepsilon=1$ applies if $k=n-1$ \cite{MSTV-96}.
Of course, $\AR(n,K_{1,k})=\lar(n,K_{1,k})$ holds
 because stars do not contain $2K_2$.
Here we prove that the presence of $2K_2$ in $\brmk$ does not
 change the value of $\lar$ if $n$ is not too small.
The particular case of $k=3$ will be reconsidered in
 Section~\ref{ss:claw+leaf}, where exact formulas for all $n$
 will be proved for further \arr\ functions.

\btm   \label{t:brmk}
For every\/ $k\geq 3$ there is a threshold\/ $n_0=n_0(k)$ such that
 $\lar(n,\brmk)=\AR(n,K_{1,k})$ holds for all\/ $n\geq n_0(k)$.
\etm

\bpf
The lower bound $\lar(n,\brmk)\geq\AR(n,K_{1,k})$ is clear
 because $K_{1,k}\subset\brmk$.
For the upper bound $\lar(n,\brmk)\leq\AR(n,K_{1,k})$ let us
 suppose $n\gg k$, and let $\psi$ be an
 edge coloring of $K_n$ with at least $\AR(n,K_{1,k})$ colors.
Then a rainbow star $S\cong K_{1,k}$ occurs, by definition;
 say vertex $v$ is its center and $u_1,\dots,u_k$ are
 its leaves, and $\psi(u_iv)=i$ for $i=1,\dots,k$.
If $S$ cannot be extended to a properly colored $\brmk$,
 we must have $\psi(u_iw)=i$ for all $w\in V(K_n)\smin V(S)$
 and all $1\leq i\leq k$.

If $\psi(xy)$ is distinct from $1,\dots,k$ for some
 $x,y\notin\{u_1,\dots,u_k\}$ (but $x=v$ or $y=v$ is allowed),
 then we find a properly colored $\brmk$ with center $x$,
 2-length leg $xyu_1$, and pendant edges $u_ix$ for $i=2,3,\dots,k$.
Otherwise the number of colors is at most $\binom{k}{2}$ on the
 edges $u_iu_j$, plus $k$ on the remaining edges.
But $\binom{k}{2}+k$ is smaller than $\AR(n,K_{1,k})$ if
 $n$ is not too small, hence the theorem follows.
\epf

We note that the condition $k\geq 3$ in the above theorem
 is necessary because
 for $k=2$ we have $\brmk\cong P_4$ with $\lar(n,P_4)=3$
 while $K_{1,2}\cong P_3$ with $\AR(n,P_3)=2$
 (see Table \ref{tab:rsmall}).

\subsection{Graphs with degrees all odd or all even}
\label{ss:deg-odd-even}

\btm {\bf\itshape (Odd Graphs)}   \label{t:odd}
\
Let\/ $G$ be an odd graph. Then
 \begin{itemize}
  \item[$(i)$] $\oar(n,G) = 1$;
  \item[$(ii)$] $\soar(n,G) = \parr(n,G) = \lpar(n,G)$
   either tends to infinity with\/ $n$, or is equal to\/ $1$
   for all\/ $n\geq n_0(G)$ sufficiently large.
 \end{itemize}
\etm

\bpf
The first part is immediately seen, as the sum of any
 even integers is even, while all vertices have odd degrees.
Concerning the second part $\soar(n,G) \geq \parr(n,G) \geq
 \lpar(n,G)$ holds by definition, and in fact the three values
 are equal whenever $G$ is an odd graph, as observed in
 Proposition~\ref{p:parity}.
We are going to prove that $\soar(n,G)=1$ holds if $n$
 is large, unless $\lim_{n\to\infty} \soar(n,G)= \infty$.

If $\soar(n,G)=1$ for some $n\geq |G|$, then $\soar(n',G)=1$
 holds for all $n'\geq n$ because a required  copy of $G$ is
 present in every $n$-vertex subgraph of $K_{n'}$ in any
 edge coloring.
Else, suppose for a contradiction that $k:=\liminf \,\soar(n,G)$
 is finite and $k>1$.
We choose $n$ with $\soar(n_0,G)=k$ and $n \geq \max_{1\leq t
 <k} R(G,t)$, the Ramsey number of $G$ with $t$ colors.
Consider any edge coloring of $K_n$.
If at least $k$ colors are used, then a strong-odd-colored copy
 of $G$ occurs, as $\soar(n,G)=k$.
If a smaller number $t$ of colors is used, then $K_n$ contains
 a monochromatic---consequently strong-odd-colored---copy
 of $G$, as $n\geq R(G,t)$.
Thus, $\soar(n,G)=1$, contradicting the assumption $k>1$.
\epf

It follows from Proposition \ref{p:num-even} that every even
 graph $G$ has $\oar(n,G) \geq \binom{|G|-1}{2}+2$.
We next show that this constant lower bound can be improved
 to linear in $n$.
Recall from Section \ref{ss:lex} that $\Lexx(n,h)=n+h(h-3)/2$
 is the number of colors in a \klek{$h$}\ coloring.

\btm {\bf\itshape (Even Graphs)}   \label{t:even}
If\/ $G$ is an even graph on\/ $k$ vertices, then\/
 $\soar(n,G) \geq \oar(n,G) \geq \Lexx(n, k-1)+1 \geq n$,
  and\/ $\cfar(n,G) \geq \Lexx(n, k-1)+1 \geq n$.
\etm

\bpf
Let $\psi$ be the \klek{$(k-1)$}\ coloring of $K_n$ with
 vertex order $v_1,\dots,v_n$.
Consider any copy of $G$, and let $v_j$ be the vertex of
 largest index in $G$.
Since $\psi(v_iv_j)=j-k+1$ holds for all
 $v_i\in V(G)\smin \{v_j\}$, and the degree of $v_j$ is even,
 $G^\psi$ does not satisfy the requirement of odd coloring,
 neither of conflict-free coloring.
Hence $\soar(n,G) \geq \oar(n,G) \geq \Lexx(n,k-1) +1 \geq n$,  as well as
 $\cfar(n,G) \geq \Lexx(n, k-1)+1 \geq n$. 
\epf

\btm {\bf\itshape (Corona of Even Graphs)}   \label{t:cor-even}
Let\/ $H$ be any even graph, and\/ $G$ obtained from\/ $H$ by
 adding a leaf to every vertex of\/ $H$.
Then for\/ $n \geq n_0(G)$, $\soar(n,G) = 1$.
\etm

\bpf
We apply Theorem \ref{t:ERCT} for $K=K_{3m}$, where
 $m = |G| = 2|H|$. 
If there is a monochromatic or rainbow copy of $K$,
 we are done, as a copy of $G$ is strong odd-colored.
Hence it suffices to consider a \Lex-colored $K$, say with
 vertices $v_1,\dots,v_{3m}$ in this order under \Lex.
(In fact, from now on it would be enough to take $3m/2$
 vertices only.)
Embed $H$ in $v_{m +1},\dots,v_{2m}$ arbitrarily.

At any $v_i$, the edges going to neighbors of higher indices
 in $H$ have mutually distinct colors.
A vertex can have either an odd or an even number of neighbors
 with lower indices, inducing a monochromatic star of odd or
 even degree, respectively.

If $v_i$ has an even number of lower neighbors in $H$, embed
 its leaf as $v_{i-m}$; and if it has an odd number
 of lower neighbors, embed its leaf as $v_{i+m}$.
In this way a strong odd coloring of $G$ is obtained.
So, no matter how many colors we use, $\soar(n, G) = 1$ holds
 for $n \geq n_0(G)$.
\epf

\subsection{Odd majority orientation and \womo}

Here we introduce a certain type of permutation on the vertex
 set, and present some consequences in connection with the
 \arr\ functions under consideration.

\bdf
Given a graph\/ $G$ on\/ $k$ vertices\/ $v_1,\dots,v_k$
 and a permutation\/ $\pi$
 (of the\/ $k!$ possible permutations) of its vertices,
 $\pi$ is called\/ \emph{odd majority orientation} if:
 \begin{enumerate}
  \item each edge\/ $e = (v_{\pi(i)}, v_{\pi(j)})$, with\/
   $\pi(i) < \pi(j)$ is oriented from\/ $v_{\pi(j)}$ to\/ $v_{\pi(i)}$;
  \item each vertex\/ $v$ has either no outgoing edges
   (i.e., $\deg^+(v) = 0$) or has an odd number of
   outgoing edges (\/$\deg^+(v) \equiv 1 \ (\mathrm{mod} \ 2)$). 
 \end{enumerate}
\edf

We also consider the following related notion.

\bdf
Given a graph\/ $G$ on\/ $k$ vertices\/ $v_1,\dots,v_k$,
 a permutation\/ $\pi$
 of its vertices is called\/ \emph{\womo} if each vertex\/
  $v_{\pi(i)}$ satisfies one of the following conditions:
 \begin{enumerate}
  \item $v_{\pi(i)}$ has either zero or an odd number of
   neighbors\/ $v_{\pi(j)}$ with\/ $\pi(j) < \pi(i)$;
  \item $v_{\pi(i)}$ has an even number of
   neighbors\/ $v_{\pi(j)}$ with\/ $\pi(j) < \pi(i)$, and no
   neighbors\/ $v_{\pi(j')}$ with\/ $\pi(i) < \pi(j')$. 
 \end{enumerate}
\edf

By definition, every odd majority orientation $\pi$ is an
 \womo, but not vice versa.
Also, if $G$ is an even graph, then it admits no odd majority
 orientation because the vertex of highest index assigned by any
 permutation violates the condition.
But many even graphs (the cycles, for instance) admit an
 \womo.

\bex
\
\begin{itemize}
 \item[$(i)$]
For\/ $k \geq 3$, the complete graph\/ $K_k$ does not admit an
 odd majority orientation.
Indeed, in any permutation, the third vertex has exactly two
 neighbors whose indices are smaller.
 \item[$(ii)$]
If\/ $G$ is not an odd graph but has an odd majority orientation
 (e.g., $G=K_4-e$), then inserting a new veretex and joining it
 to the odd-degree vertices of\/ $G$ we obtain an even graph
 that admits an odd-even ordering.
(The highest index can be assigned to the new vertex.)
\end{itemize}
\eex

The role of odd graphs and odd majority orientations is
 explored in the following result.

\btm   \label{t:odmaj}
\
\begin{itemize}
 \item[$(i)$] 
If\/ $G$ is an odd graph, and\/ $G$ has an odd majority orientation, then\/
 $\soar(n,G) = \oar(n,G) = \parr(n,G) = \cpar(n,G) =
  \lpar(n,G) = 1$ for all\/ $n \geq n_0(G)$.
 \item[$(ii)$] 
If a graph\/ $G$ admits an \womo, then\/
 $\lpar(n,G)=1$ holds for all\/ $n \geq n_0(G)$.
 \item[$(iii)$] 
If\/ $G$ does not admit any odd majority orientation, then,
 for all\/ $n\geq |G|$,
 $\soar(n,G) \geq \parr(n,G) \geq \cpar(n,G) \geq n$, and if
 in addition\/ $G$ has no \womo,
 then also\/ $\lpar(n,G) \geq n$. 
\end{itemize}
\etm

\bpf
\
$(i)$\quad 
Since $\soar(n,G)$ dominates all the other parameters,
 it suffices to consider strong odd colorings.
Assume that $G$ has $p$ vertices, and apply
 Theorem \ref{t:ERCT} for $K_p$. 
Then for every sufficiently large $n \geq n_0(p)$, no matter
 how many colors are used in an edge coloring of $K_n$,
 there is a copy $K$ of $K_p$ whose coloring is either
 monochromatic, or rainbow, or \Lex.
This coloring pattern is inherited for any embedding of $G$
 into $K$.
Since all vertex degrees are odd, a monochromatic $G$ is
 strong-odd-colored.
Also, every rainbow graph is strong-odd-colored.
Finally, if $K$ is \Lex-colored, we choose a permutation
 $\pi$ that generates an odd majority orientation. 
Embed $G$ into $K$ in accordance with $\pi$.
Any color $j$ from the corresponding vertex $v_{j+1}$ to its
 smaller-index neighbors in $G$ occurs on an odd number of
 vertices; and the edges from $v_{j+1}$ to its higher-index
  neighbors have mutually distinct colors.
Thus, a strong odd coloring is obtained, independently of the
 number of colors in $K_n$, whenever $n$ is sufficiently large.
Consequently, $\soar(n,G) = 1$. 

\msk

\nin
$(ii)$\quad We again apply Theorem \ref{t:ERCT}.
If $n$ is sufficiently large, then any edge coloring contains a
 monochromatic or rainbow or \Lex-colored $K_p$, $p=|G|$.
The first two patterns immediately yield local-parity-colored
 copies of $G$.
If a \Lex-colored $K_p$ is found, we take an \womo\ on $V(G)$.
The monochromatic star towards the lower-indexed neighbors of a
 vertex either has odd degree, which is of the same parity as
 the non-repeated colors to the higher-indexed neighbors, or
 has even degree equal to the degree of the vertex in $G$.
Both cases are compatible with local parity coloring.

\msk

\nin
$(iii)$\quad To show that $n-1$ colors do not guarantee a
 class-parity-colored copy of $G$, nor a local-parity-colored
 copy if $G$ satisfies the extra condition, we apply \Lex.
Any copy of $G$ in $K_n$ must have a vertex, say $v_{j+1}$, with
 an even number of neighbors with lower indices, hence in the
 corresponding orientation induced by $\pi$ it has an even
 number of incident color-$j$ edges dictated by \Lex\ in $K_n$.
This star is not allowed in class parity coloring.
Moreover, if $G$ has no \womo, then any embedding of $G$ in the
 \Lex-colored $K_n$ contains a vertex $v$ where a monochromatic
 star of even degree occurs towards its lower-indexed neighbors,
 and also it has at least one higher-indexed neighbor, to which
 the color of the edge is not repeated at $v$.
This is not allowed in local-parity-coloring.
\epf

This approach is easily applicable to bipartite graphs.

\btm   \label{p:omo-bip}
Let\/ $G$ be a bipartite graph. Then

\begin{itemize}
 \item[$(i)$] 
$G$ admits an \womo, hence\/ $\lpar(n,G)=1$
 for all\/ $n \geq n_0(G)$.
 \item[$(ii)$] 
If\/ $G$ is an odd graph, and also if all even-degree vertices
 of\/ $G$ belong to the same vertex class,
  then\/ $G$ has an odd-majority orientation,
  hence\/ $\soar(n,G) = \oar(n,G) = \parr(n,G) = \cpar(n,G) =
   \lpar(n,G) = 1$ holds for all\/ $n \geq n_0(G)$.
\end{itemize}
\etm

\bpf
Let $A \cup B$ be the bipartition of $G$.
Take a permutation that enumerates first (with smallest indices)
 all the vertices from $A$, and then all the vertices from $B$. 
The vertices of $B$ only have lower-index neighbors, and those
 of $A$ are adjacent only to higher-index neighbors, respectively.
Hence, an \womo\ is obtained.
And if all vertices of even degree are in $A$, then it is also
 an odd-majority orientation.
Thus, the results of Theorem \ref{t:odmaj} apply to $G$.
\epf

Cycles, no matter if even or odd, do not admit an odd-majority
 orientation, as the vertex of highest index in any permutation
 has exactly two lower neighbors.
On the other hand, the exclusion of cycles suffices:

\bpn   \label{p:omo-tree}
Every tree and forest admits an odd-majority orientation.
\epn

\bpf
We apply induction on the number of vertices.
Let $T$ be a tree or a forest.
The assertion is trivial if $T$ has no edges.
Otherwise let $uv$ be a pendant edge, where $v$ is a leaf of $T$.
By the induction hypothesis, $T-v$ admits an odd-majority orientation.
Insert $v$ as the last vertex of $T$ with highest index, and for
 the rest of the vertices keep the odd-majority orientation of $T-v$.
Then the number of lower neighbors of $u$ remains odd, and of
 course $v$ has just one lower neighbor (hence an odd number).
\epf

It is important to note that there exist graphs admitting an
 odd-majority orientation and still having both $\cpar(n,G)$
 and $\oar(n,G)$ tend to infinity with $n$.
In this sense $\lpar(n,G)$ substantially differs from all the
 seven other \arr\ functions.

\section{Conflict-free \arr\ number is linear}
\label{s:CF-linear}

The main result of this section is a general linear upper bound
 on all $\cfar(n,G)$.
Since every conflict-free coloring is also a (weak) odd coloring,
 the theorem implies the linearity of $\oar(n,G)$ as well.

We begin with a simple observation.

\bpn   \label{p:cf-Kp}
For every\/ $p \geq $ we have\/
$\cfar(p,K_p) = \binom{p}{2} - \left\lfloor \frac{p}{2} \right\rfloor + 1
  = \left\lceil \frac{p^2}{2} \right\rceil - p + 1$.
\epn

\bpf
to obtain a non-conflict-free coloring of $K_p$, it is
 necessary and sufficient that at least one vertex $v$ has
 all its incident colors occur at least twice.
If the vertex degree $p-1$ is even, then the loss compared
 to the number $\binom{p}{2}$ of colors in a rainbow $K_p$ is
 $(p-1)/2$, and if the degree is odd, then the loss is $p/2$.
Hence, $\binom{p}{2} - \left\lfloor \frac{p}{2} \right\rfloor$
 colors do not gauarantee a conflict-free coloring, but
 more colors do.
\epf

On general graphs, a universal upper bound can be guaranteed as follows.

\btm   \label{t:CF-G_p}
Let\/ $G=(V,E)$ be a graph with\/ $p$ non-isolated vertices.
Then\/
 $$
   \cfar(n,G)\leq(p-2)(n-p)+\left\lceil \frac{p^2}{2} \right\rceil - p + 1
    =(p-2)n-\left\lfloor \frac{p^2}{2} \right\rfloor + p + 1 \,.
 $$
\etm

\bpf
We apply induction on $n$.
The anchor with $n=p$ is settled in Proposition \ref{p:cf-Kp}.

In the induction step we prove $\cfar(n+1,G)\leq\cfar(n,G)+p-2$.
To do this, consider $K=K_{n+1}$ and any of its edge colorings
 $\psi$ with at least $\cfar(n,G)+p-2$ colors.
If there is a vertex $v\in V(K)$ such that $\psi$ uses at least
 $\cfar(n,G)$ colors in $K-v$, then we are done
 by the induction hypothesis.
Otherwise each $v$ is the center of at least $p-1$
 color classes that are stars.
Let $Q_v\subset N(v)$ be a set composed by selecting one vertex
 from each star color class centered at $v$.
(If an edge $e=vw$ is a singleton color class, then $w\in Q_v$
 and $v\in Q_w$ are unique choices representing that color
 at the two ends of $e$.)

Before designing a procedure how a conflict-free copy of $G$
 is found in $K^\psi$, we make some preparations in $G$.
Let $S$ be any non-extendable independent set in $G$, and
 denote $Z:=V\smin S$.
Then every $z\in Z$ has at least one neighbor in $S$.
Assume that $X=\{x_1,\dots,x_k\}\subset S$ is a minimal
 set dominating $Z$.
If the set $W:=S\smin X$ is nonempty, let further
 $Y=\{y_1,\dots,y_l\}\subset Z$
 be a minimal set dominating $W$.
The definition of $Y$ is meaningful because $S$ is maximal and
 $G$ has no isolates, hence every $w\in W$ has a neighbor
 that must belong to $Z$ as $S$ is an independent set.
Observe that, by the minimality of $X$ and $Y$, each $x_i$ has
 a neighbor in $Z$ whose unique neighbor in $X$ is $x_i$; and
 each $y_i$ has a neighbor in $W$ whose unique neighbor
 in $Y$ is $y_i$.

For convenience we label the vertices in $X$ and $Y$ in a way
 that the degrees of the $x_i$ are non-increasing,
 $d_G(x_1)\geq d_G(x_2)\geq \cdots \geq d_G(x_k)$;
 and the degrees of the $y_i$ towards $W$ are non-increasing,
 $d_W(y_1)\geq d_W(y_2)\geq \cdots \geq d_W(y_l)$.

Next, still in $G$, we specify vertex subsets $X_i\subset Z$
 ($i=1,\dots,k$) and $Y_i\subset W$ ($i=1,\dots,l$)
 sequentially in the order of increasing subscript as follows.
Artificially setting $X_0=Y_0=\es$,
 let $X_i$ be the set of all those vertices of $Z$ that have
 no neighbors in $\bigcup_{j=0}^{i-1} X_j$; and
 let $Y_i$ be the set of all those vertices of $W$ that have
 no neighbors in $\bigcup_{j=0}^{i-1} Y_j$.

Now we are in a position to design an injective mapping
 $\eta:V\to V(K)$ that embeds $G$ into $K$ and yields a
 conflict-free-colored subgraph $H\cong G$.
First, let $v_1=\eta(x_1)$ be any vertex of $K$, and let
 $\eta(X_1)$ be an arbitrarily chosen subset of $Q_{v_1}$
 (having size $|Q_{v_1}|=|X_1|$, of course).
After that, for $i=2,\dots,k$ in this order, we select
 $\eta(X_i)$ as an $|X_i|$-element subset of
  $Q_{v_i}\!\smin \bigcup_{j=1}^{i-1} \eta(X_j)$.
At this point each $z\in\eta(X_i)$ is incident with a single
 edge of color $\psi(z\eta(x_i))$, and this color occurs
 only once at $\eta(x_i)$ as well.
Note further that those colors do not appear inside
 $V(K)\smin\eta(X)$, thus the corresponding edges remain
 single representatives of their colors even when we add
 any further vertices to $\eta(X\cup Z)$.

We complete the construction by applying a similar procedure
 for the vertices of $\eta(Y)$.
To simplify notation, let us denote $U:=V(K)\smin\eta(X\cup Z)$.
Let $\eta(Y_1)$ be any $|\eta(Y_1)|$-element subset of
 $U\cap Q_{\eta(y_1)}$; and then, for $i=2,\dots,l$ in that order,
 select $\eta(Y_i)$ as a $|Y_i|$-element subset of
  $(U\cap Q_{\eta(y_i)})\smin \bigcup_{j=1}^{i-1} \eta(Y_j)$.
Now each $w\in\eta(Y_i)$ is incident with a single
 edge of color $\psi(w\eta(y_i))$, hence ensuring that the
 conflict-free requirement is satisfied at $w$.

Since $|Q_v|\geq p-1$ holds for every $v\in V(K)$, all the
 above selections are possible, and a conflict-free coloring
 of a subgraph isomorphic to $G$ is found in $K$.
\epf

Since $\oar(n,G) \leq \cfar(n,G)$ holds for all graphs $G$ and
 all $n\geq |G|$, we also obtain:

\bcr
We have\/ $\oar(n,G)=O(n)$ for every graph\/ $G$.
\ecr

\section{Paths, cycles, and small graphs}
\label{s:maxdeg-2}
\label{s:more-small}

In this section we mostly deal with the exhaustive list of all
 graphs of at most four vertices or edges, with $P_6$ as a
  slight extension, complementing the works done in
 \cite{BGR-15,GR-16} concerning $\AR(n,G)$ where $G$ is either
 a small graph or has only small components.
Our results are then summarized in Table \ref{tab:rsmall}
 for the convenience of the reader.
Moreover, some general estimates for paths and cycles of any
 length will also be presented.

Before turning to particular types of graphs, let us introduce
 a general concept that will be useful in several proofs,
 notably in Sections \ref{ss:paths} and \ref{ss:K4-e}.

\paragraph{Locally critical colors.}

Let $\psi$ be any edge coloring of $K_n$, with any number $k$ of colors.
Call a color $i$ critical at a vertex $v$ if all edges of
 color $i$ are incident with $v$.
There are two kinds of critical color classes: a single edge
 with its private color, and a star with at least two edges.
Denote the number of the former and the latter by $s$ and $t$,
 respectively.
A single edge is critical at its both ends, a star is
 critical at its center but not at its leaves.
Hence the total number of 
 incidences is equal to
\begin{center}
\qquad
\hfill
  \#(vertex, incident critical color) = $2s+t$.
  \hfill (\ref{s:more-small}.0)
\end{center}
The relevance of critical colors becomes apparent in proofs
 by induction:

\bob   \label{ob:crit}
In an inductive proof of\/ $\phi(n,G|\cF)\leq an+b$ for an
 \arr-related function\/ $\phi$ under consideration,
 assuming that\/ $\phi(n-1,G|\cF)\leq a(n-1)+b$ has been proved,
 one may restrict attention to edge colorings\/ $\psi$ of\/ $K_n$
 such that every vertex is incident with at least\/
 $a+1$ critical colors.
Due to equality\/ $(\ref{s:more-small}.0)$,
 in this case we have
  $$
    2s+t \geq (a+1)n \,.
  $$
\eob

\paragraph{Local parity coloring.}

Let us recall that $\lpar(n,G)=1$ if $\Delta(G)\leq 2$
 (Proposition \ref{p:parity} (4)) and $\lpar(n,G)=1$ holds
 for $n$ large enough whenever $G$ is bipartite or,
  more generally, admits an \womo\ (Theorems
 \ref{t:odmaj}~$(ii)$ and \ref{p:omo-bip} $(i)$), respectively).
From the graphs considered in this section, $K_4$ is the only
 one to which none of these principles can be applied.
In fact $K_4$ is an odd graph and its $\lpar$ has a
 substantially different behavior, satisfying
 $\lpar(n,K_4)=\soar(n,K_4)=\parr(n,K_4)$ by Proposition \ref{p:parity} (2).
For these reasons, we will not discuss on $\lpar(n,G)$ for the
 various graphs $G$ below separately.

\subsection{Lower bound for paths and cycles}

Let us recall first that the \arr\ numbers of cycles have been
 asymptotically determined by Montellano-Ballesteros and Neumann-Lara
 \cite{MBNL-05} as $\AR(n,C_k) =  {(k-2)/2 +1/(k-1)}n +O(1)$;
 and the problem of $\AR(n,P_k)$ for paths has been solved
  recently by Yuan \cite{Y-21+}.
We note, however, that the coloring suggested by \erd, Simonovits
 and \sos\ as a lower bound for \arr\ $C_k$ is not applicable for
  $\lar(n,C_k)$.
This fact is worth mentioning because most of the considered
 functions defined in terms of local conditions are equal on any
 graph of maximum degree 2, but $\AR(n,G)$ and $\lpar(n,G)$
 usually have a different behavior.

For the local versions concerning these graphs we have the
 following general lower bounds.

\btm   \label{t:path-cyc-LB}
Recalling that\/ $s(n,h) = hn - \binom{h+1}{2} +1$,
\begin{itemize}
 \item[$(i)$] for\/ $k \geq 4$ we have\/ $\lar(n,C_k) = \soar(n, C_k)  =  \oar(n,C_k) = \cfar(n,C_k) \geq 
  s(n, \lfloor{(k-1)/3}) \rfloor ) +1$;
 \item[$(ii)$] for\/ $k \geq 6$ we have\/ $\lar(n,P_k) = \soar(n, P_k)  =  \oar(n,P_k) = \cfar(n,P_k) \geq 
 \parr(n,P_k) \geq \cpar(n,P_k) \geq 
  s(n, \lfloor{(k-3)/3} \rfloor ) +1$.
\end{itemize}
\etm

\bpf
Recall that for any graph $G$ with maximum degree at most 2
 the four functions $\lar$, $\soar$, $\oar$, $\cfar$ are equal.
Moreover, the inequalities $\soar(n,G)\geq \parr(n,G)\geq
 \cpar(n,G)$ are valid for every graph $G$.
Therefore, we only have to prove lower bounds on $\lar(n,C_k)$
 and $\cpar(n,P_k)$.
The constructions for the two cases are quite similar, but there
 are some differences in the details.
However, in either case, the idea is that the constructed edge
 coloring does not contain any paths and cycles above a certain
 length, hence it suffices to restrict attention to the smallest
 length relevant for a formula.

\msk

\nin
$(i)$\quad
Let $k \geq 4$, $k \equiv 1$ (mod 3), and consider the
 \krs{(k-1)/3}-coloring of $K_n$.
In any properly colored cycle $C$ of length $\ell$, at most two
 consecutive vertices can occur in the monochromatic part $B$,
 otherwise there would be a vertex of the
 cycle with two incident edges in $B$ having the same color.
Hence, if $\ell \geq 3|A|+1$, then $C$ is not properly colored.

\msk

\nin
$(ii)$\quad
Let $k \geq 6$, $k \equiv 0$ (mod 3), and consider the
 \krs{(k-3)/3}-coloring of $K_n$.
Also here, in any properly colored path $P$ of length $\ell$,
 at most two consecutive vertices can occur in the monochromatic
 part $B$, otherwise there would be a vertex of the
 path with two incident edges in $B$ having the same color.
Hence, if $\ell \geq 3|A|+3$, then the edges of $P$ inside $B$
 form a linear forest with at least one component of length
 exceeding 1.
This is not allowed in class parity coloring.

\msk

The inequalities proved on $\ell$ verify the lower bounds in
 both parts of the theorem.
\epf

\subsection{The cycle $C_4$}
\label{ss:C4}

\btm   \label{t:C4}
We have\/
$\lar(n,C_4) = \soar(n,C_4) = \oar(n, C_4) = \cfar(n,C_4) =
 \parr(n,C_4) = \cpar(n,C_4) = n +1 \,.$
\etm

\bpf 
Due to Observation \ref{ob:ineq} that the first four values are
 equal on any cycle; moreover they provide an upper bound on the
 last two values for any graph.

The lower bound on $\lar(n,C_4)$ is the particular case $k=4$ of
 Theorem \ref{t:path-cyc-LB}, using the \krs{1} coloring.
Although it does not work for the functions involving parity
 conditions, here an alternative coloring can be defined:
In fact the same number of colors is achieved with the
 substantially different \klek{3}, as well.
In \klek{3} the vertex of highest index determines a color class
 $P_3$ with exactly two edges of the same color in any copy
 of $C_4$, hence no class-parity-colored $C_4$ occurs.

For the upper bound, let $\psi$ be any edge coloring of $K_n$
 with more than $n$ colors.
Selecting one edge from each of the first $n+1$ color classes
 we obtain a rainbow graph with more edges than vertices.
The same inequality also holds in at least one connected
 component of this selection.
Such a component contains more than one cycle,
 hence we can select a connected rainbow subgraph
 whose structure is one of the following:
 \begin{itemize}
  \item[$(a)$] two vertices connected by three paths $P,P',P''$,
   any two of which are internally vertex-disjoint;
  \item[$(b)$] two vertex-disjoint cycles $C',C''$
   connected by a path $P$;
  \item[$(c)$] two cycles $C',C''$ sharing precisely one
   vertex.
 \end{itemize}
As a matter of fact, we can reduce $(a)$ to the following
 favorable case:
 \begin{itemize}
  \item[$(d)$] some rainbow cycle $C$ of even length contains
   no repeated color.
 \end{itemize}
It is clear that $(a)$ reduces to $(d)$ because any two of
 $P,P',P''$ form a rainbow cycle, and the lengths of (at least)
 two of those paths have the same parity.

If $(d)$ holds, we assume that $C$ is a shortest rainbow
 cycle of even length.
If $C$ is a 4-cycle, then we are done.
Otherwise let $u,v$ be two vertices at distance 3 along $C$.
They are connected by two paths along $C$, say $P'=uxyv$ and
 $P''=C\smin\{x,y\}$.
Now consider the edge $e=uv$ of $K_n$.
If $\psi(uv)\neq\psi(ux)$ and also $\psi(uv)\neq\psi(yv)$,
 then an odd-colored $C_4$ is found on $\{u,x,y,v\}$.
Otherwise $\psi(uv)$ is absent from $P''$, thus $P''\cup e$
 is a rainbow even cycle shorter than $C$, a contradiction.

In case $(b)$ we assume that $P$ is as short as possible.
Let $P$ have its ends $u\in V(C')$ and $v\in V(C'')$.
Consider the path $uxyv$ where $ux$ is an edge of $C'$ and
 $yv$ is an edge of $C''$.
If $\psi(xy)$ does not occur in $P$, then omitting the
 (at most one) edge of color $\psi(xy)$ from $P\cup C'\cup C''$
 and inserting the edge $xy$ we see that a rainbow subgraph
 of type $(a)$ can be found, and the proof is done.
Otherwise, if $\psi(xy)$ occurs in $P$, the minimality
 condition on $P$ implies that $P$ is the single edge $uv$
 (as $xy$ alone would also connect $C'$ with $C''$),
 and since $\psi(xy)=\psi(uv)$ does not occur in $C'\cup C''$,
 we obtain an odd-colored $C_4$ on $\{u,x,y,v\}$.

It remains to analyze $(c)$ where both $C'$ and $C''$ are
 odd rainbow cycles.
Now we assume that $|C'|+|C''|$ is as small as possible.
If $|C'|>3$ (i.e., at least~5), let $uxyv$ be a subpath of $C'$
 disjoint from $V(C'')$.
As in the proof for $(d)$, we can consider the color $\psi(uv)$
 of edge $e=uv$ and either find that there is a rainbow $C_4$
 on $\{u,x,y,v\}$ or obtain the contradiction that $C'$
 can be shortened to $C'\smin\{x,y\}$ by the insertion of $e$.
The same argument applies to $C''$, as well.
As a consequence, $C'\cup C''$ is the bow-tie graph $K_1+2K_2$.
Let its two 3-cycles be $wux$ and $wvy$; by asumpton all their
 six edges have mutually distinct colors.
To complete the proof, we consider the edge $uv$ of $K_n$.
If $\psi(uv)\notin\{\psi(uw),\psi(yv)\}$, then there is an
 odd-colored $C_4$ on $\{u,w,y,v\}$.
Likewise, if $\psi(uv)\notin\{\psi(ux),\psi(wv)\}$, then
 there is an odd-colored $C_4$ on $\{u,x,w,v\}$.
But at least one of these two cases must hold because
 $C'\cup C''$ is a rainbow graph, implying
 $\{\psi(uw),\psi(yv)\}\cap\{\psi(ux),\psi(wv)\}=\es$.

Thus, an odd-colored $C_4$ is found in every $\psi$.
\epf

\subsection{Short paths}
\label{ss:paths}

For $n\geq 5$ we have $\AR(n,P_4)=3$, as proved in \cite{BGR-15}.
Three colors are not sufficient for $n=4$, as shown by the
 proper edge 3-coloring of $K_4$.
However, three colors suffice for all the other parameters
 considered in this paper.

\bpn
We have $\soar(n,P_4) = \oar(n,P_4) = \cfar(n,P_4) = \lar(n,P_4)
 = \parr(n,P_4) = \cpar(n,P_4) = 3$
  for all\/ $n\geq 4$.
\epn

\bpf
A monochromatic spanning star in color 1 with a monochromatic
 $K_{n-1}$ in color 2 shows for all but $\lpar(n,P_4)$ that
 two colors are not sufficient.
Due to $\AR(n,P_4)=3$ for $n\geq 5$, we only have to verify the
 tightness of the lower bound 3 for $n=4$.
If at least three colors are used in $K_4$, consider the two
 edges $v_1v_2$ and $v_3v_4$.
They are connected by an edge whose color is distinct from both
 $\psi(v_1v_2)$ and $\psi(v_3v_4)$, hence a properly colored
 $P_4$ is found.
\epf

We also have tight results for paths of length four and five,
 as follows.

\btm
We have\/
 $\lar(n,P_5)=\oapf=\soapf=\cfapf=\parr(n,P_5)=\cpar(n,P_5)=4$
 for every\/ $n\geq 5$.
\etm

\bpf
We have already seen that $\lar(n,P_5)=\oapf=\soapf=\cfapf$
 and $\lar(n,P_5)\geq \parr(n,P_5)\geq \cpar(n,P_5)$.
So, it will suffice to prove $\cpar(n,P_5)>3$ and
 $\lar(n,P_5)\leq 4$.

For the lower bound, let $\psi$ be the edge 3-coloring of
 $K_n$ where all edges incident with $v_1$ have color 1,
 all edges incident with $v_2$ except $v_1v_2$ have color~2,
 and all edges not meeting $\{v_1,v_2\}$ have color 3.
Consider any copy $P$ of $P_5$.
One end of $P$ should be $v_1$, otherwise $P$ contains a
 monochromatic $P_3$ (with $v_1$ in its middle), not allowed
 in class parity coloring, and the proof is done.
If $v_2$ is not the other end of $P$, then a monochromatic $P_3$
 with middle $v_2$ occurs.
If the two ends of $P$ are $v_1$ and $v_2$, then the internal
 three vertices induce a monochromatic $P_3$ in
 $K_n^\psi-v_1-v_2$.
Hence, $K_n^\psi$ does not contain any class-parity-colored $P_5$.

For the upper bound, let $n\geq 5$, and assume that $\psi$ is
 an edge coloring of $K_n$ without any proper $P_5$.
We begin with two simple observations.
 \begin{itemize}
  \item[$(a)$] If $P$ is a proper $P_4$ with an end vertex $u$
   whose incident edge in $P$ has color $i$, then all vertices
   $w\notin V(P)$ have $\psi(uw)=i$.
  \item[$(b)$] $K_n^\psi$ does not contain any proper $C_4$.
 \end{itemize}
Here $(a)$ simply expresses that $P$ has no extension to a
 proper $P_5$.
To see $(b)$, let $C=pqru$ be a proper $C_4$, and consider
 any vertex $w\notin V(C)$.
Since both $P'=pqru$ and $P''=rqpu$ are proper $P_4$, by $(a)$
 we should have $\psi(pu)=\psi(uw)=\psi(ru)$, contradicting
 the assumption that $C$ is properly colored at~$u$.

Assume now that $\psi$ uses at least four colors.
Since $\AR(P_4)=3$, we know that a rainbow $P_4$ occurs; say,
 the path $P:=vxyz$ has color pattern $(1,2,3)$ (meaning
  $\psi(vx)=1$, $\psi(xy)=2$, $\psi(yz)=3$).
In the sequel we write $w$ (sometimes with a subscript) for
 vertices not contained in $P$.

Due to $(b)$ we can assume $\psi(vz)=1$.
(More precisely $(b)$ implies color 1 or 3 on $vz$, but
 1 can be taken by symmetry reasons.)
Then, applying $(a)$ for the proper paths $zyxv$, $vzyx$, and $vxyz$, we obtain
 \begin{itemize}
  \item $\psi(vw)=1$, $\psi(xw)=2$, $\psi(zw)$ for all $w\notin V(P)$.
 \end{itemize}

We consider the possible positions of an edge of color 4.
 \begin{itemize}
  \item If $\psi(vy)=4$, then $wxvyz$ is a proper $P_5$ with
   color pattern $(2,1,4,3)$ for any $w$.
  \item  If $\psi(xz)=4$, then $(a)$ implies $\psi(yw)=3$ for
   all $w$ by the path $vxzy$;
    and then $wyxzv$ is a proper $P_5$ with
   color pattern $(3,2,4,1)$ for any $w$.
  \item If $\psi(yw)=4$ for some $w$, then $vxywz$ is a proper
   $P_5$ with color pattern $(1,2,4,3)$.
  \item If $n\geq 6$ and $\psi(w_1w_2)=4$, then $vxw_1w_2z$ is
   a proper $P_5$ with color pattern $(1,2,4,3)$.
 \end{itemize}
Hence, in all possible cases a proper (in fact, rainbow) $P_5$
 has been found.
\epf

\brm
There is a substantial difference between the behaviors of\/
 $P_4$ and\/ $P_5$, as\/ $\Ar{P_4}=3$ is a constant, while\/
 $\Ar{P_5}=n+1$ is linear in\/ $n$, despite that the other
 three functions for\/ $P_5$ remain constant\/ $4$.
Although in the analyzed four positions of an edge of color\/ $4$
 we always found a rainbow\/ $P_5$, this seeming contradiction
 arises because\/ $(a)$ has been applied, which is not valid
 under the requirements of\/ $\AR$.

Also, a substantial difference occurs between\/ $P_5$ and\/ $P_6$,
 demonstrated by Theorem \ref{t:P6} below, as the values
 are constant for the former grow linearly for the latter.
\erm

\btm   \label{t:P6}
For every\/ $n\geq 6$ we have
 $\oar(n,P_6) = \soar(n,P_6) = \cfar(n,P_6) = \lar(n,P_6) =
  \parr(n,P_6) = \cpar(n,P_6) = n+1$.
\etm

\bpf
We know that all the first four functions are equal, and that
 they provide an upper bound on the last two.
The lower bound $n+1$ is verified by the \krs{1} coloring,
 namely a rainbow spanning star with monochromatic $K_{n-1}$
 on its leaf set.
Any copy of $P_6$ contains a subpath with more than one edge from
 the monochromatic $K_{n-1}$, hence not class-parity-colored.

The upper bound is more complicated to prove.

\msk

\nin
\underline{The start:}
 \ $n=6$, \ $\lar(6,P_6)=7$.

\ssk

Let $\psi$ be a coloring of $K_6$ with at least 7 colors.
Since $\AR(6,P_4\cup P_2)=7$ by Proposition 6.3 of \cite{BGR-15},
 we can label the vertices in such a way that $P=vxyz$ is a
 $P_4$ with color pattern $(1,2,3)$ and $uw$ is an edge of
 color 4 that joins the other two vertices.
Trying to avoid a proper $P_6$ under $\psi$, step by
 step we obtain restrictions on the colors of edges in $K_6$,
 which eventually will force the presence of a proper $P_6$ anyway.

First we show that
 \begin{itemize}
  \item $\psi(vz)=1$
 \end{itemize}
  can be assumed (or 3, but the two cases are symmetric).
Indeed, if $\psi(vz)\notin\{1,3\}$ then $vxyz$ is a proper
 $C_4$, so e.g.\ viewing the two 4-cycles $yzvx$ and $vzyx$
 we would obtain $\psi(xu)\in\{1,4\}\cap\{2,4\}$, hence it
 should be color 4.
Repeating this for all vertices of the 4-cycle it would follow
 that all edges incident with $u$ and $w$ should have color 4,
 but then $uvxyzw$ would be a proper $P_6$.

Next, in a similar way, the three paths $zyxv,vzyx,vxyz$ imply
 \begin{itemize}
  \item $\psi(vu),\psi(vw)\in\{1,4\}$,
  \item $\psi(xu),\psi(xw)\in\{2,4\}$,
  \item $\psi(zu),\psi(zw)\in\{3,4\}$.
 \end{itemize}
So far only four colors have been used, hence three new colors
 must occur, but only four edges are uncolored, namely
 $xz,yv,yu,yw$; only one of the four can have an old color
  from $\{1,2,3,4\}$.
There must occur an old color in both $wuyxzv$ and $uwyxzv$,
 hence the edge of the old color is
  $\{uy,xz\}\cap\{wy,xz\}=xz$ with color $\psi(xz)\in\{1,2\}$.
The other three colors are new and mutually distinct:
 \begin{itemize}
  \item $\psi(yv)=5$,
  \item $\psi(yu)=6$,
  \item $\psi(yw)=7$.
 \end{itemize}

Now the path $vwyzxu$ takes its color pattern from the alternatives
 of $(\{1,4\},7,3,\{1,2\},\{2,4\})$ which is a non-proper $P_6$
 only if
 \begin{itemize}
  \item $\psi(xu)=\psi(xz)=2$.
 \end{itemize}
But then $vxuwyz$ is a rainbow $P_6$ with color pattern
 $(1,2,4,7,3)$.

\msk

\nin
\underline{The induction step:}
 \ $n>6$.

\ssk

Consider now $n \geq 7$, assuming that $\lar(n-1,P_6)=n$
 has been proved.
Let $\psi$ be any edge coloring of $K_n$ with some $k>n$
 colors, and assume that no~properly colored $P_6$ is present.
Denoting by $s$ and $t$ the number of critical edge- and
 critical star-classes, respectively, based on Observation
 \ref{ob:crit} we may assume $2s+t\geq 2n$, or $s+t/2\geq n$.

If $s+t \geq n+2$, then we build a rainbow subgraph $F$ of $K_n$
 using one edge from each of the first $n+2$ critical color classes.
The sum of degrees in $F$ is equal to $2(n +2)$, hence
 the average degree is $2(n+2)/n < 3$ as $n \geq 7>4$.
Consequently there is a vertex in $F$ with degree at most 2,
 and deleting it from $K_n$ we obtain $K_{n -1}$ with at least
 $n$ colors, hence containing a properly colored $P_6$.

From now on we can assume $s+t \leq n+1$.
Together with $s+t/2 \geq n$ this yields $s\geq n$ or $s=n-1$,
 the latter implying $t = 2$ also.

Assume first $s \geq n$.
We then consider a rainbow $F$ with $n$ critical edges.
If there is a $P_4$ in $F$, then both ends of this $P_4$ can be
 extended, to obtain a properly colored $P_6$. 
If there is no $P_4$ in $F$, then we apply the fact that
 $\ex(n, P_4) \leq n$, and equality holds if and only if
  $3\mid n$ and $F\cong\frac{n}{3}K_3$.
So, for $n \equiv 1, 2$ (mod~3) we are done. 
For $n \equiv 0$ (mod 3) and $n\geq 6$ there are two
 vertex-disjoint triangles formed by critical edges.      
Take an edge connecting them.
It must be of a distinct color, and a rainbow $P_6$ is obtained.

Finally it remains to settle $s=n-1$ and $t = 2$.
We concentrate on the graph $F$ formed by the $n-1$
 critical edges.
If $P_4\subset F$, and also if $2K_3\subset F$, the proof is
 done as above.
Similarly, if $3P_2\subset F$, then the three disjoint
 critical edges can be joined by two further edges to form a
 properly colored $P_6$.
In particular, it follows that if no proper $P_6$ are found,
 then $F$ either is connected or has exactly two components.

If $F$ is connected but contains no $P_4$, then
 $F\cong K_{1,n-1}$, say with center $v$, and $K_n-v$ is
 colored with at least two colors, all distinct from the
 star centered at $v$.
Inside $K_n-v$ we take a 2-colored $P_3=xyz$ and any two
 further vertices $u,w$.
Then $xyzvuw$ is a properly colored $P_6$.

If $F$ has two components, then one of them is $K_3$ (two trees
 would have only $n-2$ edges), the other one is $K_{1,n-4}$.
Then any edge connecting the triangle with a leaf of the star
 creates a properly colored $P_6$.
\epf

\subsection{Claw with leaf}
\label{ss:claw+leaf}

\def \clap {K_{1,3}^+}

Let us introduce the notation $\clap$ for the graph obtained by
 attaching a pendant edge to one leaf of $K_{1,3}$.

\bpn
$\lar(n,\clap)=\AR(n,\clap)=\lfloor n/2 \rfloor + 2$
 for all $n\geq 6$.
\epn

\bpf
It is proved in \cite{BGR-15} that $\lfloor n/2 \rfloor + 2$
 is an upper bound on $\AR(n,\clap)$ whenever $n\geq 6$, hence
  also for $\lar(n,\clap)$.
To see that $\lar(n,\clap)$ is not smaller, we take a rainbow
 matching of size $\lfloor n/2 \rfloor$ and assign a new color
 $c$ to all the other edges of $K_n$.
Then in any copy of $\clap$ the vertex of degree~2 is incident
 with at least two edges of color $c$, hence
 $\lfloor n/2 \rfloor + 1$ colors are not enough to guarantee
 a properly colored $\clap$.
\epf

\btm
$\soar(n,\clap) = \oar(n,\clap) = \parr(n,\clap) =
 \cpar(n,\clap) = 2$ for\/ $n \geq 7$. 
\etm

\bpf
Clearly, the monochromatic $K_n$ does not satisfy the conditions
 on $\clap$ for any of the four functions, hence 2 is a valid
 lower bound.

Since $\soar(n,\clap)$ dominates the other three functions, the
 proof will be done if we show the upper bound $\soar(n,\clap) \leq 2$.
Let the vertices of $K_n$ be $v_1,\dots,v_n$, where $n\geq 7$.
Suppose for a contradiction that $\psi$ is a non-monochromatic
 edge coloring without strong $\clap$.

Suppose first there is a vertex with at least three edges of the
 same color at a vertex, say
  $\psi(v_1v_2) = \psi(v_1v_3) = \psi(v_1v_4) = 1$.
Then, by the exclusion of a strong $\clap$, all edges
 between $\{v_2,v_3,v_4\}$ and $\{v_5,\dots,v_n\}$ have color 1.
But then $v_2,v_3,v_4$ have degree at least $n-3\geq 4$ in
 color 1, therefore they also are adjacent to each other in color 1.
Hence they have degree $n-1$ in this color, which implies by a
 similar reason that the entire $K_n$ is monochromatic.

Otherwise every color occurs on at most two edges at each vertex.
Then, since $n\geq 7$, at least three colors occur at $v_n$.
Say $\psi(v_1v_n)=1$, $\psi(v_2v_n)=2$, $\psi(v_3v_n)=3$.
Viewing $v_n$ as the possible center of a $\clap$ we obtain that
 $\psi(v_1v_i)=1$ must hold for all $4\leq i\leq n$.
Hence the degree of $v_1$ in color 1 is high, and we are back
 to the case which has already been settled.
\epf

\brm
Let us note that for\/ $\soar(n,\clap)<3$ the condition\/
 $n\geq 7$ is necessary.
This is shown by the edge coloring of\/ $K_6$ where color\/ $1$
 induces\/ $K_{3,3}$ and color\/ $2$ induces\/ $2K_3$.
On the other hand, e.g.\ $\soar(n,\clap)<3$ is easy to prove
 for all\/ $n\geq 5$.
A good start is to take a\/ $2$-colored\/ $P_3$, say\/ $v_4v_1v_5$,
 then reduce the problem to\/ $n=5$ by picking any\/ $v_2,v_3$,
 observe that the triangle\/ $v_1v_2v_3$ should be monochromatic
 unless a strong-odd-colored\/ $\clap$ occurs,
 get a further\/ $2$-colored\/ $P_3$, then another
 monochromatic triangle, and so on, until a required copy of\/
 $\clap$ is found.
\erm

\btm
$\cfar(n,\clap)= 3$.
\etm

\bpf
Assign color 1 to all edges incident with a selected vertex $v$
 of $K_n$, and color 2 to all edges of $K_n-v$.
Then no Cf-colored $\clap$ occurs, proving $\cfar(n,\clap)\geq 3$.

On the other hand, if an edge coloring of $K_n$ uses at least
 three colors, we can find a Cf-colored $P_4=wxyz$.
Supplementing it with an edge $vx$, where $v$ is a fifth vertex,
 we obtain a Cf-colored $\clap$.
\epf

\subsection{Triangle with leaf, $K_1+(K_2\cup K_1)$}
\label{ss:K3+leaf}

In order to simplify notation, let us introduce $\khp$ for the
 graph often called ``paw'' or ``pan'' in the literature,
 obtained from $K_3$ by attaching a pendant edge.
We begin with a formula that can easily be deduced from known
 earlier results.

\btm
$\lar(n,\khp)=n$.
\etm

\bpf
Since every triangle has to be rainbow not only under the
 $\AR(n,G)$ scenario but also under $\lar(n,G)$, we obtain
  $$
    n = \AR(n,K_3) \leq \lar(n,\khp) \leq \AR(n,\khp) = n \,, 
  $$
 hence equality holds throughout.
\epf

\btm   \label{t:oar-khp}
We have\/ \hbox{$\soar(n,\khp) = \oar(n,\khp) =
 \parr(n,\khp) = \cpar(n,\khp) = 2$} for all\/ $n\geq 6$.
For smaller\/ $n$ we have\/ $\soar(4,\khp) = \soar(5,\khp) = 4$
 and\/ $\oar(4,\khp) = \oar(5,\khp) = 2$.
\etm

\bpf
One color is not enough, since $\khp$ contains vertices of
 degree 2 (against $\oar$ and $\soar$) and also vertices of
 opposite parity (against $\cpar$ and $\parr$).
Also, if $n=5$, assigning color 1 and color 2 equally to the
 four edges incident with a selected vertex and color 3 to the
 other six edges shows that three colors are not enough to
 guarantee a strong-odd-colored $\khp$.
If $n=4$, we omit one vertex incident with color 3 from the
 5-vertex construction.
Hence the values of $\soar(n,\khp)$ and $\oar(n,\khp)$ cannot be
 smaller than claimed.

For upper bounds, note that $\soar$ dominates all the three
 other functions, hence it will suffice to find a
 strong-odd-colored copy of $\khp$.
We first settle the cases where some vertex is
 incident with at least three edges of the same color.
Let $v$ be a vertex where the number of
 monochromatic edges is largest.
Say, $v_0$ has $d\geq 3$ neighbors $v_1,\dots,v_d$
 adjacent to $v_0$ in color 1.
Then $\psi(v_iv_j)=1$ must hold for all $1\leq i<j\leq d$,
 otherwise the triangle $v_0v_iv_j$ supplemented with any
 further edge $v_0v_k$ is a strong-odd-colored $\khp$.
This yields a $K\cong K_{d+1}$ in color 1.
For $n=4$ no more colors would be possible.
Note further that, for any large $n$, there are no edges
 of color 1 between $V(K)$ and $V(K_n)\smin V(K)$, since $d$ has
  been chosen to be the largest degree among all color classes.

Consider a vertex $x\notin V(K)$, and let $\psi(v_0x)=2$.
If two further vertices $v_i,v_j$ have $\psi(v_ix)=\psi(v_jx)=2$
 then the triangle $v_iv_jx$ with the pendant edge $v_0x$ forms
 a desired $\khp$.
The same conclusion holds if two vertices $v_i,v_j$ have
 $\psi(v_ix)=3$ and $\psi(v_jx)=4$.
Thus, it follows that $d+1=4$, and we have
 $\psi(v_0x)=\psi(v_1x)=2$ and $\psi(v_2x)=\psi(v_3x)=3$.

In particular, if $n=5$, then the assumption $d\geq 3$ restricts
 the number of colors to 3.
If $n\geq 6$, we take a sixth vertex $y\notin V(K)\cup\{x\}$.
If $\psi(xy)$ is 2 or~3, then the triangle $v_0v_1x$ or
 $v_2v_3x$ with the pendant edge $xy$ forms a desired $\khp$;
 and if $xy$ has another color (1 or 4), then we can take the
 triangle $v_1v_2x$ for the same.

Hence, from now on we can assume that each color has degree
 at most~2 at every vertex.
For $n\geq 6$ this yields at least $\lceil (n-1)/2 \rceil\geq 3$
 colors at each vertex.
If a 3-colored $K_3$ occurs, we supplement it with a pendant
 edge whose color is distinct from the colors of the two
 incident edges of that triangle at the attaching vertex.
If none of the triangles has three colors, but some color occurs
 twice at some vertex, we consider a monochromatic $P_3$.
Say, $\psi(wx)=\psi(wy)=1$, and a further edge $yz$ has color~2.
Repeatedly applying the ``degree at most two'' and ``no rainbow
 triangle'' conditions, we obtain:
 $\psi(wz)=2$ (for $w$ and $wyz$), $\psi(x,z)=1$ (for $z$ and
 $wxz$);
 and choosing a fifth vertex $u$ such that $\psi(ux)=3$,
 we derive $\psi(ux)=3$ (for $x$ and $uwx$).
And then the color of $uz$ cannot be defined, because it should
 be 2 or 3 for the triangle $uwz$, but already color 2 occurs
 twice at $z$ and color 3 occurs twice at $u$.
Thus, a 3-colored triangle is unavoidable, and the proof is
 done for $n\geq 6$.

To see $\soar(4,\khp)\leq 4$, consider any coloring
 of $K_4$ with four or more colors, and select one edge
 from each of the first four colors.
If the missing two edges form $P_3$, then a rainbow $\khp$ has
 been obtained.
If they form $2K_2$, then we have a rainbow $C_4$.
Insert one of its diagonals, and omit the edge of the same color
 from $C_4$ if it is present there, or otherwise omit any edge
 of $C_4$.
Then again a rainbow $\khp$ is obtained.

Finally, to see $\soar(5,\khp)\leq 4$, consider any coloring
 of $K_5$ with four colors $1,2,3,4$ or more.
If the removal of a vertex $v$ does not destroy any of $1,2,3,4$,
 then $K_5-v\cong K_4$ is colored with at least four colors and
 a rainbow $\khp$ can~be found as above.
Otherwise two of the five vertices destroy the same color,
 say~4, which is then a single edge.
Moreover, we know that each color has degree at most 2 at each vertex.
So, if the removal of $v_1,v_2,v_3$ destroys color $1,2,3$,
 respectively, then colors $1,2,3,4$ occur on at most seven
 edges altogether.
Hence there are at least five colors, and a 3-colored triangle
 $C_3$ occurs.
From the union of the two color classes which do not appear in
 this $C_3$, at least one edge meets $C_3$, thus yielding
 a rainbow $\khp$ and completing the proof.
\epf

\btm   \label{t:cf-khp}
$\cfar(n,\khp) = 3$.
\etm

\bpf
To show that two colors not suffice, assign all but one edges of
 $K_n$ with color 1, and one edge colored 2.   
Then at least one vertex of $K_3$ in any copy of $\khp$ will be
 incident with a monochromatic star, not conflict-free-colored.

Assume that an edge coloring $\psi$ of $K_n$ uses at least
 three colors.
If $n\geq 5$, we apply the fact that $\AR(n,P_4)=3$ holds
 for all $n\geq 5$.
Consider a rainbow $P_4 = wxyz$ where $\psi(wx) = a$,
 $\psi(xy) = b$, $\psi(yz) =  c$.
A Cf-colored $\khp$ is immediately found unless $\psi(wy)=a$
 and $\psi(xz)=c$.
Otherwise, if $\psi(wy)=a$ and $\psi(xz)=c$, the color $\psi(wz)$
 differs from at least one of $a$ and $c$; say, it is not $a$.
Then the edges $wx,wy,wz,xy$ induce a Cf-colored $\khp$.

Finally, let $n=4$.
If $\psi$ uses more than three colors, we find a rainbow triangle
 (on applying the fact $\AR(n,C_3)=n$)
 and supplement it with any pendant edge.
If just three colors are used, pick one edge from each color class.
If they form a triangle or a $P_4$, we are done as above.
If they form the star $K_{1,3}$ with center $w$ and leaves
 $x,y,z$, then consider $\psi(xy)$.
If it is $c$, we have found a Cf-colored
 (in fact $\lar$-colored) $\khp$.
And if it is $a$ (or $b$), then $xywz$ (or $yxwz$) is a rainbow
 $P_4$, and we are done as above.
\epf

\subsection{The diamond $K_4-e$}
\label{ss:K4-e}

\btm
$\oar(n, K_4 -e )  = 3$.
\etm

\bpf
An edge coloring of $K_n$ with all but one edges colored 1 and one edge colored 2 
shows that at least three colors are needed.

Suppose an edge coloring $\psi$ uses at least three colors.

First we show that if a rainbow triangle $C_3=xyz$ with colors
 $\psi(xy)= a$, $\psi(vz) = b$, $\psi(zx) = c$ occurs, then
 an odd-colored $K_4-e$ can be found.
Consider a fourth vertex $w$ outside $C_3$.
If it is adjacent to $x ,y$  with distinct colors, we are done.  
So we may assume $\psi(wx) = \psi( wy)=d$, and also
 $\psi(wx) = \psi(wz)=d$ for the same reason.
If $d$ is not one of $a,b,c$ then delete any edge from the
 rainbow $C_3$ and we are done.  
Otherwise assume without loss of generality that $d = a$.
Delete the edge colored $a$ from the rainbow $C_3$ and we are done.  

Next, if $n\geq 5$, similarly to the proof of Theorem \ref{t:cf-khp}
 we consider a rainbow $P_4 = wxyz$ where
 $\psi(wx) = a$, $\psi(xy) = b$, $\psi(yz) =  c$.
By the above, we may
  assume that there is no rainbow $C_3$ under $\psi$.
Then $\psi(wy)$ is $a$ or $b$, and $\psi(xz)$ is $b$ or $c$.
But then, independently of $\psi(wz)$, omitting the edge $xy$
 we obtain an odd-colored $K_4-e$ on $\{w,x,y,z\}$.

Hence, we are left with the smallest case $n=4$.
We next observe that a 3-colored $K_{1,3}$ also yields an
 odd-colored $K_4-e$.
Indeed, assume $\psi(wx)=a$, $\psi(wy)=b$, $\psi(wz)=c$.
Excluding a rainbow $C_3$, we may assume $\psi(xy)=a$.
Then an odd-colored $K_4-e$ is found unless $\psi(xz)=c$.
But then inserting the edge $yz$ and omitting $wx$ yields an
 odd-colored $K_4-e$, independently of $\psi(yz)$.

Consider now any $\psi$ on the edges of $K_4$ with at least
 three colors.
Assume $\psi(wx)=a$, $\psi(wy)=b$, and consider a third color $c$.
There are three possible positions of a color-$c$ edge:
 $xy$ or $wz$ or and edge from $z$ to $\{x,y\}$.
This yields a rainbow $C_3$, or $K_{1,3}$, or $P_4$, respectively.
The proof is complete.
\epf

\btm
For every\/ $n\geq 4$ we have\/ $\oard=\parr(n,K_4-e)=n+1$.
\etm

\bpf
Since $\soar\geq\parr$ is universally valid, we need to prove
 $\parr(n,K_4-e)>n$ and $\soar(n,K_4-e)\leq n+1$.
The lower bound is provided by the \krs{1} coloring:
 a rainbow spanning star and monochromatic $K_{n-1}$ in a new color.
Then every copy of $K_4-e$ has at least three vertices in the
 monochromatic part, and any induced subgraph of $K_4-e$ with
 more than two vertices contains a vertex of degree 2.
But a 2-regular color class is not allowed in a strong parity
 coloring of $K_4-e$ because at the degree-3 vertices some
 color class must have an odd degree.

The proof of the upper bound is by induction on $n$.
The basic case of $\sod{4}=5$ is obvious because using five colors
 on the edges of $K_4$ only one repetition occurs, and removing
 one edge from the duplicated color we obtain a rainbow $K_4-e$.

Consider now $n \geq 5$, assuming that $\soar(n-1,K_4-e)=n$
 has been proved.
Let $\psi$ be any edge coloring of $K=K_n$ with some $k>n$
 colors, and assume that no strong $K_4-e$ is present.
Denoting by $s$ and $t$ the number of critical edge- and
 critical star-classes, respectively, based on
 Observation \ref{ob:crit} we may assume $2s+t\geq 2n$.

We are going to prove that the number of single-edge color classes
 is at most $2n/3$.
More explicitly, each connected component in the graph formed
 by the single-edge colors is $K_2$ or $P_3$.
For this, we need to exclude $P_4$ and $K_3$ from the graph
 of critical single edges.
The exclusion of a $P_4=vxyz$ is immediate because
 by the single-edge criticality of the edges in $P_4$ we have
  $\psi(vy),\psi(xz)\notin\{\psi(vx),\psi(xy),\psi(yz)\}$,
 hence $P_4$ would yield a strong $K_4-e$.
In case of a $K_3=xyz$ with three critical edges, if
 $\psi(xw)\neq\psi(yw)$ for a $w\notin\{x,y,z\}$, then
 a strong $K_4-e$ would occur on $\{w,x,y,z\}$.
Consequently we should have $\psi(xw)=\psi(yw)=\psi(zw)$; but
 then $K_3-e$ together with $w$ would form a strong $K_4-e$.

As a consequence, we have $s\leq \lfloor 2n/3 \rfloor$,
 implying $s+t\geq \lceil 4n/3 \rceil := n^*$.
This also implies $n\geq 6$, because for $n=5$ we would have
 at most 3 critical edges and then the number of edges
 should be at least $s+2t = 2(s+t) - s \geq 4n-3s \geq 11$,
 while $K_5$ has just 10 edges.

For $n\geq 6$ we consider $n^*$ color classes of $\psi$
 (any choice) and select one edge from each of them.
In this way a rainbow graph, say $F$, with degree sum
 $2n^* < 3n$ is obtained.
Consequently, $\delta(F)\leq 2$.
Let $v$ be a vertex of minimum degree in $F$.
Then $F-v$ is a rainbow graph with at least $n^*-2\geq n$
 edges.
Thus, $\psi$ has more colors in $K-v=K_{n-1}$ than the number
 of vertices, and a proper $K_4-e$ must occur by the
 induction hypothesis.
This final contradiction completes the proof of the theorem.
\epf

\btm   \label{t:l-K4-e}
For every\/ $n\geq 4$ we have
 $$\lfloor (3n-3)/2 \rfloor + 1 \leq \lard \leq 2n-3.$$
\etm

\bpf
\
\textit{Lower bound.}\quad
The construction
 is a slight modification of the \Lex\ or the \klek{3} coloring,
 depending on the parity of $n$.

We can artificially say that $\lod{2}=2$ and $\lod{3}=4$,
 because the trivial 1-coloring of $K_2$ and the rainbow $K_3$
 exhibit the largest possible numbers of colors that we can
 use for $n=2$ and $n=3$ without a proper $K_4-e$.

Having a coloring of $K_{n-2}$ at hand, say on the vertex set
 $Z$, we adjoin two new vertices $x,y$ and construct a coloring
 for $K_n$ on $\{x,y\}\cup Z$.
We use three new colors: one for the edge $xy$, one for the
 star of $x$--$Z$ edges, and one for the
 star of $y$--$Z$ edges.

A proper subgraph $H\cong K_4-e$ should contain at least one
 vertex outside $Z$ and at least two vertices inside $Z$.
But all $H$ in such a position contain at least two $x$--$Z$
 edges or at least two $y$--$Z$ edges, hence not properly colored.

\msk

\nin
\textit{Upper bound.}\quad
We procered by induction on $n$.
The basic case of $\lod{4}=5$ is obvious because using five colors
 on the edges of $K_4$ only one repetition occurs, and removing
 one edge from the duplicated color we obtain a rainbow $K_4-e$.
(This is the same as the basic case for strong odd coloring.)

Let $\psi$ be any edge coloring of $K=K_n$ with some $k\geq 2n-3$
 colors, and assume that no strong $K_4-e$ is present.
Recall that a color $i$ is critical at a vertex $v$ if
 all edges of color $i$ are incident with $v$.
As discussed at the beginning of this section, a
 critical color class is either a single edge or a star.

If a vertex $v$ is incident with at most two critical colors,
 then $\psi$ uses at least $2n-5$ colors in $K-v=K_{n-1}$,
  hence a strong $K_4-e$ occurs by the induction
 hypothesis, contradicting the assumptions.
It follows that each vertex is incident with
 at least three critical colors.

Construct a mixed graph $H=(V,E,A)$ as follows.
 \begin{itemize}
  \item $V=V(K)$.
  \item $vw\in E$ if $vw$ is a single-edge class in $\psi$.
  \item $\overrightarrow{vw}\in A$ if $\psi(vw)$ is a critical
   color at $v$, and the color class is a star centered at $v$
    with more than one edge.
 \end{itemize}
For each vertex $v$ we denote by $N_i(v)$ the set of vertices $x$
 such that $vx\in E$ or $\overrightarrow{vx}\in A$, the
 neighbors for which $v$ is critical in color $i$, and
 set $\displaystyle N^*(v):=\{v\}\cup\bigcup_i N_i(v)$
  for their union together with $v$ itself.

Consider any $v$.
Say the critical colors at $v$ are $1,\dots,k$.
The edges of the complete $k$-partite graph
 $\langle N_1(v),\dots,N_k(v) \rangle$ are monochromatic,
 otherwise  it would contain two incident edges whose other
 ends are in distinct classes $N_i(v)$, hence a 2-colored
 $P_3=xyz$ would occur and with $v$ it would form a rainbow
 $K_4-e$, contradicting the assumptions.
This one ``crossing'' color is not critical because $k>2$.
We denote this color by $c(v)$.

As a consequence, there can be three types of critical colors
 at a vertex $x\in N_i(v)$, namely
  \begin{itemize}
   \item[$(a)$] $vx$ is an edge critical for both $v$ and $x$;
   \item[$(b)$] $x$ is the end or center of a critical edge or
    star entirely inside $N_i(v)$;
   \item[$(c)$] there is a color $p>k$ and a vertex $y\in N_p(x)$
    such that $p\neq c(v)$ and $y\notin N_i(v)$ for any
    $1\leq i\leq k$.
  \end{itemize}

Suppose that case $(c)$ holds for some $v,x,y$.
We then select a vertex $z\in N_2(v)$.
Observe that $c(v)\notin\{1,2,p\}$ and $\psi(vy)\notin\{1,p\}$.
We have arrived at the contradiction that $\{v,x,y,z\}$
 contains a proper $K_4-e$.
Hence, the proof will be done if we show that case $(c)$ is
 unavoidable.

Consider the sets $N^*$ for all vertices of $K_n^\psi$, and let
 $v$ be a vertex such that $|N^*(v)|$ is as small as possible.
Assume that $N^*(v)\supset N_1(V)\cup N_2(v)\cup N_3(v)$.
Select now a vertex $x\in N_1(v)$.
The color $c(v)$ is not critical, therefore
 $N_2(v)\cup N_3(v)\subset N^*(v)$
 is a nonepty set disjoint from $N^*(x)$.
However, $|N^*(x)|\geq |N^*(v)|$ holds by assumption,
 implying the presence of a vertex $y\in N^*(x)\smin N^*(v)$.
This completes the proof.
\epf

The next result dealing with class parity coloring demonstrates a
 less expected application of the \erd--Rado Canonical Theorem.

\btm
$\cpar(n,K_4-e)=5$ if $n$ is large.
\etm

\bpf
The following 4-coloring shows that
 $\cpar(n,K_4-e)\geq 5$ holds for every~$n$.
Select two vertices $v_1,v_2$ in $K_n$, and denote
 $K'=K_n-v_1-v_2$.
Assign color 1 to all edges from $v_1$ to $K'$ and color 2 to
 all edges from $v_2$ to $K'$.
Let $K'$ be monochromatic in color 3, and assign color 4 to the
 edge $v_1v_2$.

Consider any copy $G$ of $K_4-e$.
If $G\subset K'$ is monochromatic, then it contains vertices of
 degree 2 and 3 in the same color class, not allowed in
 class parity coloring.
If both $v_1,v_2\in V(G)$, then at least one of colors 1 and 2
 induces a $P_3$ color class, not allowed either.
Finally, if exactly one of $v_1,v_2$ is in $G$, say $v_1$, it
 can have degree 2 or 3 in $G$.
Degree 2 yields a monochromatic $P_3$ color class as above.
Degree 3 yields the odd graph $K_{1,3}$ in color 1.
But then $G-v_1\cong P_3$ is monochromatic in color 3, not allowed.

The proof of the opposite inequality $\cpar(n,K_4-e)\leq 5$
 requires more work.
Let $\psi$ be any edge coloring of $K_n$ with at least five
 colors.
We assume that $n$ is sufficiently large to guarantee a $K_4$
 with rainbow or \Lex-colored or monochromatic coloring
 via Theorem \ref{t:ERCT}.
If this $K_4$ is rainbow, it contains a rainbow
 $K_4-e$ and we are done.
If it is \Lex-colored, we take a vertex order of $K_4-e$
 where the two vertices of degree 2 are in the middle, and
 the two degree-3 vertices have the lowest and highest index.
Then \Lex\ generates a color class $K_{1,3}$ at the highest
 vertex, and two single-edge color classes, hence a strong
 odd coloring, and we are done also in this case.

From now on we may assume that there is no rainbow $K_4$ and
 no \Lex-colored $K_4$ under $\psi$.
Let $K$ be a largest monochromatic complete subgraph in $K_n$,
 say in the highest color $k\geq 5$.
The choice of $n$ guarantees $|K|\geq 4$.

We next analyze the colors from the external vertices
 $v_i\in V(K_n)\smin V(K)$ to $K$.
There must be at least one edge of some color $c_i\neq k$
 from $v_i$ to $K$, because $K$ is not extendable to a larger
 monochromatic complete graph.
If there are two such edges of distinct colors
 $c_i\neq c_i'\neq k$, then we complete them
 to a copy of $K_4-e$ with a further vertex in $K$.
Indeed, this $K_4-e$ has a monochromatic $K_3$ (an even graph)
 and two single edges as color classes, hence a
 class parity coloring is obtained.
So, we may assume that each $v_i$ is adjacent to $K$ with edges
 either colored $k$ or colored with the same color~$c_i$
 other than $k$.
This also implies that at least three vertices
 $v_1,v_2,v_3$ exist outside $K$, because $\psi$ uses
 at least five colors.

Next, we observe that there is at most one edge of color $k$
 from $v_i$ to $K$.
Otherwise we find a $C_4$ of color $k$ with a diagonal of
 color $c_i$, hence a class-parity-colored $K_4-e$.
So, we may assume that each $v_i$ is adjacent to $K$ with
 at least $|K|-1$ edges in color $c_i$.

If $c_i=c_j$ for some $i\neq j$, then the neighborhoods of
 $v_i$ and $v_j$ in color $c_i$ share at least $|K|-2$
 vertices, that is at least two.
In this case we find a $C_4$ of color $c_i$ with a diagonal of
 color $k$, hence again a class-parity-colored $K_4-e$ occurs.
So, we may assume without loss of generality that $c_i=i$
 for $i=1,2,3$.

If $\psi(v_1v_2)=1$ (or likewise an edge of color 2 or 3 is
 assigned to an edge incident with $v_2$ or $v_3$, respectively,
 inside $\{v_1,v_2,v_3\}$) then we pick two vertices $w_1,w_2$
 from $K$ such that $\psi(v_iw_j)=i$ for all $i,j\in\{1,2\}$.
Then the sequence $(w_1,w_2,v_2,v_1)$ generates a \Lex-colored
 $K_4$, a contradiction.
And if none of the above occurs, but color 1 is present as
 $\psi(v_2v_3)=1$, then we pick a $w$ from $K$ such that
 $\psi(v_iw)=i$ for $i=1,2,3$.
(There are at least $|K|-3\geq 1$ possible choices for $w$.)
Then $\{w,v_1,v_2,v_3\}$ produces a properly colored $K_4-e$,
 as only $v_1v_2$ and $v_1v_3$ may possibly have the same color
 among incident edges.

We are left with the case where none of the three edges inside
 $\{v_1,v_2,v_3\}$ get any colors from $1,2,3$.
If at least two colors are used there, then suitably omitting
 an edge we obtain a rainbow $K_4-e$.
If the triangle is monochromatic, then we supplement it with
 two edges incident with $w$ and obtain a parity-class-colored
 copy of $K_4-e$.
\epf

\subsection{The complete graph $K_4$}
\label{ss:K4}

Our last result on conflict-free colorings is valid for both
  $K_4-e$ and $K_4$, and the two can be handled together.

\btm
We have\/
 $\cfar(n, K_4 -e) = \cfar(n,K_4) =  n+1$.
\etm

\bpf
The lower bound is provided by \klek{3} using $n$ colors.
It begins with a rainbow triangle, and each later vertex
 has a monochromatic star backward.
Hence the last vertex of any $K_4-e$ or $K_4$ does not have
 a locally singleton color.

The upper bound follows by finding a rainbow $C_4$ in any
 edge coloring that uses more than $n$ colors in $K_n$ via
 Theorem \ref{t:C4}, and extending it to a \cfcc\
  $K_4-e$ or $K_4$, on applying
 Proposition \ref{p:deg2} (iii).
\epf

Now we turn to the other functions on $K_4$.
Since $K_4$ is an odd graph, we have $\oar(n,K_4) =
 \lpar(n,K_4) = \parr(n,K_4) = 1$ for all $n$, and $\cfar(n, K_4)$
 has been determined together with $\cfar(n, K_4-e)$.
Below we determine the growth order of $\lkn$, and give a
 lower bound on $\soar(n, K_4)$.

\btm   \label{t:l-K4}
We have\/ $\lkn = \Theta(n^{3/2})$ as\/ $n\to\infty$.
More explicitly, for every\/ $n\geq 4$,
 $$ \displaystyle
   (1/2 - o(1))\, n^{3/2} = \ex(n,C_4)+2 \leq \lkn \leq
     n \left\lceil \sqrt{2n} \ \right\rceil
 $$
\etm

\bpf
For the lower bound, let $G$ be any $C_4$-free graph of order $n$.
Assign mutually distinct colors to the edges of $G$, and
 extend this to $K_n$ by assigning a fresh new color to all
 edges of $\overline{G}$.
A $K_4\subset K_n$ cannot be proper because from the new edges
 it may only contain a matching (either just one edge or a $2K_2$),
 but then a $C_4$ would be composed from edges of the
 rainbow $G$, which cannot be the case.

The upper bound is obvious if $n=4$, because then
 six colors make $K_4$ rainbow, and 6 is much smaller
 than $4\cdot\left\lceil \sqrt{2\cdot 4} \ \right\rceil$.
For larger $n$ we apply induction, assuming that the upper bound
 is valid for $n-1$.

Let $\psi$ be any coloring of $K_n$ with
 $f(n):=n \left\lceil \sqrt{2n} \ \right\rceil$ colors.
Then
 $$ \textstyle
 f(n)/n = \left\lceil \sqrt{2n} \ \right\rceil \leq
 \left\lceil \sqrt{2n} \ \right\rceil +
 (n-1)(\left\lceil \sqrt{2n} \ \right\rceil
 -\left\lceil \sqrt{2n-2} \ \right\rceil) = f(n) - f(n-1)$$
 holds for every $n$.
Consequently, we can assume that there are at least
 $\sqrt{2n}+1$ critical colors each vertex, for otherwise
 an obvious induction applies.

We now define a digraph $D=(V,A)$ with vertex set $V=V(K_n)$
 and arc set $A$ in the following way.
First select one edge from each critical color.
Those edges form a rainbow undirected graph $G=(V,E)$.
Second, if an edge $e=uv\in E$ is a color class itself,
 take two arcs $\overrightarrow{uv},\overrightarrow{vu}\in A$
 for $e$.
Third, if an edge $e=uv\in E$ represents a star color class,
 orient $e$ towards the center of the star, as an arc in $A$.
In this way $D$ has at least $n\sqrt{2n} + n$ arcs.
The in-degree of each vertex is at least $\sqrt{2n} + 1$;
 but we now concentrate on the out-degrees $d_1,\dots,d_n$.

We say that an unordered vertex pair $\{u,v\}$ is assigned to
 a vertex $w$ if both $\overrightarrow{wu},\overrightarrow{wv}\in A$.
The number of vertex pairs assigned to the
 $i^\mathrm{th}$ vertex is equal to $\binom{d_i}{2}$.
Since $\sum_{i=1}^n d_i=|A|\geq n\sqrt{2n} + n$,
 applying Jensen's inequality we obtain
  $$
    \sum_{i=1}^n \binom{d_i}{2}
     \geq n \cdot \binom{|A|/n}{2}
      \geq n \cdot \frac{(\sqrt{2n} + 1)\sqrt{2n}}{2}
       > n^2 > 2 \binom{n}{2}.
  $$
Thus, there exist a pair $\{u,v\}$ assigned to at least three
 vertices, say $x,y,z$.

The six edges between $\{u,v\}$ and $\{x,y,z\}$ have
 mutually distinct colors as they originate from the rainbow
 graph $G$; and each of their color classes has
 its center in $\{u,v\}$, hence none of them occurs on the
 three edges inside $\{x,y,z\}$.
The color $\psi(uv)$ may occur as one of those six colors
 connecting $\{u,v\}$ with $\{x,y,z\}$, say its one end is $z$.
Then $\{u,v,x,y\}$ induces a $K_4$, which is either rainbow
 or 5-colored where the only one
 color coincidence is $\psi(uv)=\psi(xy)$.
In either case a proper $K_4$ is found.
\epf

Currently we do not have a strong upper bound on
 $\soar(n, K_4)$; the lower bound is linear in $n$, while
 the upper bound grows with $n^{3/2}$.

\btm   \label{t:s-K4}
For every\/ $n\geq 4$ we have\/ $\lar(n,K_4) \geq \soar(n, K_4)
 = \parr(n, K_4) = \lpar(n, K_4) \geq 2n-2$, and also\/
 $\parr(n, K_4) \geq \cpar(n, K_4) \geq 2n-2$.
\etm

\bpf
We have seen that the upper bound $\lar(n,K_4)$ and also the
 inequalities are valid by the hierarchy of criteria defining
 the corresponding \arr\ functions.
Moreover, the claimed equalities follow from
 Proposition \ref{p:parity} (2) as $K_4$ is an odd graph.

For the lower bound $2n-2$, consider a coloring of $K_n$
 with a rainbow star at $v_n$, and take a \Lex\ coloring
 of $K:=K_n-v_n$.
The total number of colors used is $2n -3$. 
If a copy of $K_4$ contains the center of the rainbow star,
 then the vertex of largest index under \Lex\ in $K$ has exactly
 two edges of the same color (forming a monochromatic $P_3$)
 and an edge of another color to $v_n$, hence this $K_4$
 is neither local-parity-colored nor class-parity-colored.
Otherwise the copy is a $K_4$ under \Lex\ on $n-1$ vertices,
 and the vertex of second highest index violates the condition.
\epf


\renewcommand{\arraystretch}{1.2}
\begin{table}[!ht]
\hspace*{-2.7cm}
\footnotesize
\begin{tabular}{cccccccc}
\hline
$G$ & $\AR(n, G)$ & $\lar(n,G)$ & $\soar(n,G)$ & $\cfar(n,G)$ & $\oar(n,G)$ & $\parr(n,G)$ & $\cpar(n,G)$ \\
\hline
\hline
$P_2$ & 1 & 1 & $\equiv$ & $\equiv$ & $\equiv$ & $=$ & $=$ \\
\hline
$P_3$ & 2 & 2 & $\equiv$ & $\equiv$ & $\equiv$ & $=$ & $=$ \\
$2P_2$ & 2 & 1 & $\equiv$ & $\equiv$ & $\equiv$ & $=$ & $=$ \\
\hline
$P_4$ & 3 & 3 & $\equiv$ & $\equiv$ & $\equiv$ & $=$ & $=$ \\
$P_3 \cup P_2$ & 3 & 2 & $\equiv$ & $\equiv$ & $\equiv$ & $=$ & $=$ \\
$K_{1,3}$ & $\lfloor \frac{n}{2} \rfloor + 2$ & $\lfloor \frac{n}{2} \rfloor + 2$ & 1 & 2 & 1 & $=$ & $=$ \\
$3P_2$ & $n + 1$ & 1 & $\equiv$ & $\equiv$ & $\equiv$ & $=$ & $=$ \\
$C_3$ & $n$ & $n$ & $\equiv$ & $\equiv$ & $\equiv$ & $=$ & $=$ \\
\hline
$K_{1,3}$ + leaf & $\lfloor \frac{n}{2} \rfloor + 2$ & $\lfloor \frac{n}{2} \rfloor + 2$ & 2 & 3 & 2 & $=$ & $=$ \\
$K_{1,3} \cup P_2$ & $\lfloor \frac{n}{2} \rfloor + 2$ & $\lfloor \frac{n}{2} \rfloor + 2$ & 1 & 2 & 1 & $=$ & $=$ \\
$K_3$ + leaf & $n$ & $n$ & 2 & 3 & 2 & $=$ & $=$ \\
$P_3 \cup 2P_2$ & $n + 1$ & 2 & $\equiv$ & $\equiv$ & $\equiv$ & $=$ & $=$ \\
$C_3 \cup P_2$ & $n + 1$ & $n$ & $\equiv$ & $\equiv$ & $\equiv$ & $=$ & $=$ \\
$P_4 \cup P_2$ & $n + 1$ & 3 & $\equiv$ & $\equiv$ & $\equiv$ & $=$ & $=$ \\
$P_5$ & $n + 1$ & 4 & $\equiv$ & $\equiv$ & $\equiv$ & $=$ & $=$ \\
$2P_3$ & $n + 1$ & $n +1$ & $\equiv$ & $\equiv$ & $\equiv$ & $=$ & $=$ \\
$K_{1,4}$ & $n + 2$ & $n+2$ & 2 & $=$ & $=$ & $=$ & $=$ \\
$C_4$ & $\lfloor \frac{4n}{3} \rfloor$ & $n +1$ & $\equiv$ & $\equiv$ & $\equiv$ & $=$ & $=$ \\
$4P_2$ & $2n - 1$ & 1 & $\equiv$ & $\equiv$ & $\equiv$ & $=$ & $=$ \\
\hline
$P_6$ & $n + 2$ & $n+1$ & $\equiv$ & $\equiv$ & $\equiv$ & $=$ & $=$ \\
$K_4-e$ & LB: $\ex(C_3, C_4) + 2$ & LB: $\lfloor \frac{3n-1}{2} \rfloor$ & $n+1$ & $n+1$ & 3 & $n+1$ & 5 \\
 & UB: $\ex(C_3, C_4) + n + 1$ & UB: $2n-3$ &&&&& \\
$K_4$ & $\lfloor \frac{n^2}{4} \rfloor +2$ & LB: $\ex(C_4)+2$ & LB: $2n-2$ & $n+1$ & 1 & $= \soar$ & LB: $2n-2$ \\
  && UB: $n \left\lceil \sqrt{2n} \ \right\rceil$ & UB: $n \left\lceil \sqrt{2n} \ \right\rceil$ &&&& UB: $\soar$ \\
\hline
\end{tabular}
\caption{\Arr\ numbers (possibly except for very small $n \leq 8$), odd and
 parity versions: exact values for all graphs
 having at most four edges and for $P_6$, and currently known best estimates
 for the other two graphs with four vertices.
LB = lower bound, UB = upper bound, $\ex(G_1, G_2)=\ex(n, \{G_1, G_2\})$;
 ``$\equiv$'' means that the equalities
 $\soar=\oar=\cfar=\lar$ hold, as Observation \ref{ob:ineq} (3) applies
 to $G$;
 ``$=$'' means that the value is equal to the preceding entry of the same row,
 but not necessarily due to a general structural principle.
The eighth function has $\lpar(n,G)=1$ for all $G\neq K_4$
 (and $\lpar(n,K_4)=\parr(n,K_4)=\soar(n,K_4)$ by Proposition \ref{p:parity}~(2)), either for all $n$ or for all
 sufficiently large $n$.
   \label{tab:rsmall}}
\end{table}


\section{Concluding remarks and open problems}
\label{s:conclude}

In this concluding section we collect many representative
 problems concerning the functions considered in this paper.
Evidently our choice is subjective, yet we hope that the
 enclosed list below opens a wide area of research with
 many more interesting results and problems to follow.

Table \ref{tab:rsmall} summarizes the main results of
 Section \ref{s:more-small} and can serve as the benchmark for
 further research and also as the start for induction steps in
 case of generalizing some graphs into wider families of graphs.

Our list of problems is organized into subsections concentrating on:
 \begin{enumerate}
  \item completion of results concerning specific graphs and small graphs;
  \item further understanding of the hierarchy of the parameters;
  \item algorithmic complexity problems concerning Odd Majority Ordering;
  \item problems concerning the effect of graph operations;
  \item general host graphs, instead of complete graphs, and hypergraph version.
 \end{enumerate}

\subsection{Problems on specific graphs}

\bpm
Determine tight asymptotics for\/ $\lar$, $\soar$, $\oar$,
  $\cfar$, $\parr$, $\cpar$
 for paths and cycles.
In particular, find reasonable lower and upper bounds on\/
 $\parr(n,C_k)$ and\/ $\cpar(n,C_k)$.
(See Theorem \ref{t:path-cyc-LB}.)
\epm

\bpm
 Find a construction for\/ $\cfar(n,K_k)$ that supplies
  a linear lower bound\/ $c_k\, n$.  
\epm

\bpm
 Determine all the parameters for\/ $K_4-e$ and\/ $K_4$ to
 complete the table of small graphs\/
  $($Table \ref{tab:rsmall}$\,)$. 
\epm

\bpm
 Find sharper---or even exact---bounds for\/
  $\soar(n, K_k)$ and\/ $\lar(n, K_k)$,
  whose order is asymptotically determined.
\epm

\bpm
Determine $\phi(|G|,G)$ for arbitrary graphs\/ $G$ and for any\/ 
 $\phi\in\{\AR,\lar,\soar,\oar,\cfar,\parr,\cpar,\lpar\}$.
\epm

\bpm
Compute the missing parameters of small graphs for small\/ $n$.
\epm

Concerning Theorem \ref{t:brmk} we raise:

\bpm
Is it true that for every\/ $q$ the spider\/ $S=S_{q*2,r*1}$
 with\/ $q$ legs\/ $P_3$ and\/ $r$ legs\/ $P_2$ (pendant edges)
 satisfies\/ $\lar(n,S)=\AR(n,K_{1,q+r})$ whenever\/ $r$ is
 sufficiently large with respect to\/ $q$\,?
\epm

\subsection{Problems on the hierarchy of constraints}

Table \ref{tab:hierarc} presents the hierarchy between the eight
 coloring requirements as implied by the definitions and proved
 in Observation \ref{ob:ineq} and Proposition \ref{p:parity}.
For instance, a rainbow coloring (AR) satisfies all the other
 conditions as well, whereas an odd coloring (OD) or a local
 parity coloring (LP) is not guaranteed to fulfill any of the
 other seven.
In particular (OD, LP) is a simple example of an incomparable pair.
Interesting questions arise both concerning comparable and
 concerning incomparable pairs.

\paragraph{Comparable classes.}

Positioned higher in the hierarchy implies that the corresponding
 parameter is not smaller.
For two types $\cF_1,\cF_2$ of coloring constraints, and a graph
 $G$, let us write $f(n,G|\cF_1) \Rightarrow f(n,G|\cF_2)$ if
 $f(n,G|\cF_1) \geq f(n, G|F_2)$ holds for all $n \geq n_0 = n_0(G)$.

\bpm
If\/ $f(n,G|\cF_1) \geq f(n, G|\cF_2)$ for all\/ graphs $G$,
 then does there exist a\/ $G'$ and\/ $n_0=n_0(G')$ such that\/
 $f(n,G|\cF_1) > f(n,G|\cF_1)$ holds for all\/ $n \geq  n_0$?
\epm

Some positive cases are immediately read out from
 Table \ref{tab:rsmall}, but not all.
Concerning $\soar(n,G) \geq \parr(n,G)$ we know of no graphs for
 which $\soar(n,G) > \parr(n,G)$ at least for all large $n$.
Such a strict inequality can hold only if $G$ is an even graph,
 as otherwise, by Proposition \ref{p:parity} (2),
 $\soar(n,G) = \parr(n,G)$.
Furthermore if $\parr(n,G)$ is realized by odd color classes,
 then still equality holds.
This fact leads to the following problem.

\bpm
Does there exists an even graph\/ $G$ and an edge coloring\/
 $\psi=\psi(n)$ of $K_n$ with\/ $\parr(n,G)$ or more colors for
 every large\/ $n$ in which a copy of\/ $G$
 occurs with even color classes but not with odd color classes? 
\epm

For such a graph $\soar(n,G) > \parr(n,G)$ would hold. 
With reference to Observation \ref{ob:ineq}~(3) let us note that
 even the following more restricted case remains open.

\bpm
Prove or disprove:
If\/ $G$ has maximum degree at most\/ $2$, then\/
 $\parr(n,G) = \soar(n,G)$, or even\/ $\cpar(n,G) = \soar(n,G)$.
\epm

We mention that for $G = K_4-e$, the values satisfy the
 inequalities $\soar(n,G) = \parr(n,G)
 = n+1 > \cpar(n,G) = 5 > \lpar(n,G) = 1$;
 these resolve the other relations. 

Also we have $\soar(n,K_4 -e) = n+1 > \oar(n,K_4-e) = 3
 > \lpar(n.K_4-e ) = 1$ while
 $\oar(n,K_4) = 1$ and $\lpar(n,K_4 ) \geq 2n-2$,
 showing that OD and LP are incomparable.

\paragraph{Incomparable classes.}

The diagram in Table \ref{tab:hierarc} indicates that the
 defining properties of the pairs
 (CF, SOD), (CF, SP), (CF,CP), (CF, LP),
 (OD, SP), (OD,CP), (OD, LP), and (CP, LP) are incomparable.
Numerically we also know that SOD, SP, CP, LP can have a
 quadratic growth, while CF and OD have a linear upper bound,
 hence a quadratic gap occurs, established e.g.\ by~$K_5$.

In the other direction, the table contains several graphs
 with $\cfar(n,G)>\soar(n,G)$.
So it remains to settle the status of (CP, LP) and of the
 pairs involving OD.
Concerning the former, there are lots of examples satisfying
 $\cpar(n,G)>\lpar(n,G)$ and $\oar(n,G)>\lpar(n,G)$,
  since $\Delta(G)\leq 2$ implies
  $\lpar(n,G)=1$ by Proposition \ref{p:parity} (4).
In fact,
 $\cpar(n,G)-\lpar(n,G)>cn$ and and $\oar(n,G)-\lpar(n,G)>cn$
 can hold for arbitrarily large $c$, establishing any large
 linear gap, as proved for long paths in
 Theorem \ref{t:path-cyc-LB} $(ii)$.
However, the following three cases remain open.

\bpm
Prove or disprove: There exist graphs\/ $G_1,G_2,G_3$ such that
 \begin{itemize}
  \item[$(i)$] $\oar(n,G_1)>\parr(n,G_1)$,
  \item[$(ii)$] $\oar(n,G_2)>\cpar(n,G_2)$,
  \item[$(iii)$] $\lpar(n,G_3)>\cpar(n,G_3)$.
 \end{itemize}
\epm

We have not even a single example of such $G_1,G_2,G_3$, neither
 a proof that such graphs do not exist.
(Certainly a $G_1$, if exists, would also serve as $G_2$.)
Clarification of these three cases would settle whether the
 hierarchy exhibited in Table \ref{tab:hierarc} coincides with
 the Hasse diagram of the partial order among the eight
 parameter classes under study.

\subsection{Problems on parity-driven vertex orders}

\bpm
Determine the complexity of the following decision problem:

\msk

{\sc Odd-Majority Orientation ( OMO ) :}

{\bf Input:} An undirected graph\/ $G=(V,E)$.

{\bf Question:} Does\/ $G$ admit an odd-majority orientation?

\msk

\nin
Also, once the answer is affirmative on\/ $G$, how much time
 does it take to find an odd-majority orientation?
How many permutations of\/ $V$ and how many orientations of\/ $E$
 correspond to odd-majority orientations of\/ $G$\,?
\epm

Membership of {\sc OMO} in \np\ is clear.
The problem is of interest not only in general but also for
 special classes of graphs.
According to Theorem~\ref{p:omo-bip} $(ii)$ the answer is
 affirmative whenever $G$ is a bipartite odd graph.
Such graphs can be recongnized in $O(|V|+|E|)$ time, linear in
 the input size.
The proof of the proposition shows that a suitable permutation
 can also be found in linear time, as it only needs to obtain
 the bipartition of $G$.

\bpm
Solve the analogous questions of algorithmic and enumerative
 nature on \womo s:

\msk

{\sc Odd-Even Ordering ( OEO ) :}

{\bf Input:} An undirected graph\/ $G=(V,E)$.

{\bf Question:} Does\/ $G$ admit an \womo?
\epm

In Theorem \ref{p:omo-bip} we observed that in a bipartite graph
 it is sufficient for the existence of an odd-majority orientation
 that all vertices of even degree belong to the same vertex class.
The following question arises as a possible natural extension.

\bpm
Let\/ $G$ be a bipartite graph, in which the set of even-degree
 vertices is independent.
Does\/ $G$ admit an odd-majority orientation?
\epm

\subsection{Problems on the effect of graph operations}

Motivated by the Adding Edge Lemma (Proposition \ref{p:deg2})
 and the observations in Section \ref{ss:jump},
 we raise the following problem.

\bpm
Concerning all the eight parameters\/ $\phi(n,G)$, where\/
 $\phi\in\{\AR,\lar,\soar,\oar,\cfar,\parr,\cpar,\lpar\}$,
 taken over all graphs\/ $G$ with\/ $|G|\leq n$,
  \begin{itemize}
   \item[$(i)$] determine the maximum of\/
    $\phi(G+e) - \phi(G)$ where\/
     $e\notin E(G)$ is any new edge, if\/ $\phi$ allows this
      difference to be positive;
   \item[$(ii)$] determine the maximum of\/
    $\phi(G) - \phi(G+e)$ where\/
     $e\notin E(G)$ is any new edge, if\/ $\phi$ allows this
      difference to be positive;
   \item[$(iii)$] the above two values under the restriction
    that\/ $e$ joins two vertices of degree\/ $2$ in\/ $G$,
  \end{itemize}
 as a function of\/ $n$.
\epm

As we have seen in Proposition \ref{p:deg2}, in some cases the
 considered parameters cannot increase; and $\AR(n,G)$ is an
 obvious example where it cannot decrease by edge insertion.
A large jump can also occur in $\cpar(n,G)$ by joining two
 vertices of degree 2, as shown by $\cpar(n,K_4-e)=5$ and
 $\cpar(n,K_4)\geq 2n-2$.

More generally we ask:

\bpm
Study the effect of further graph operations on the functions\/
 $\phi(n,G)$.
\epm

\subsection{Other host graphs and hypergraphs}

One branch of \arr\ theory deals with the so-called
 \emph{rainbow number}.
Given two graphs, $G$ and $H$, where $H$ serves as host graph,
 the goal is to determine the smallest number $k=k(G,H)$ of colors
 such that every edge $k$-coloring of $H$ contains a rainbow
 subgraph isomorphic to $G$.
Analogous problems can be raised concerning the seven functions
 introduced here.

\bpm
Given two graphs\/ $G$ and\/ $H$, and a coloring type\/
 $\Phi\in\{\mathrm{LR}$, $\mathrm{SOD}$, $\mathrm{OD}$,
  $\mathrm{CF}$, $\mathrm{SP}$, $\mathrm{CP}$, $\mathrm{LP}\}$,
 determine the smallest integer\/ $k=k(G,H;\Phi)$ such that
 every edge coloring of\/ $H$ with at least\/ $k$ colors
 contains a copy of\/ $G$ whose induced coloring is of type\/ $\Phi$.
\epm

Moreover, it is very natural to seek hypergraph analogues of
 interesting graph problems.
Recall that a hypergraph (finite set system) is $r$-uniform if
 all its edges have exactly $r$ elements.
Beyond the edge colorings of $K_n$, one may consider edge colorings
 of the complete $r$-uniform hypergraph $\mathcal{K}_n^{(r)}$,
 whose edge set consists of all $r$-element subsets of an
 $n$-element vertex set.
Given a fixed $r$-uniform hypergraph $\mathcal{H}$, and a property
 $\mathcal{P}$ of edge colorings,
 one can ask for the minimum number $\phi_{\mathcal{P}}(n,\mathcal{H})$
 of colors such that every edge coloring of $\mathcal{K}_n^{(r)}$
 with at least $\phi_{\mathcal{P}}(n,\mathcal{H})$ colors contains a copy
 of $\mathcal{H}$ that satisfies property $\mathcal{P}$.
Similarly, if $\mathcal{H}_0$ is an $r$-uniform host hypergraph
 that contains at least one copy of $\mathcal{H}$, one can define
 $\phi_{\mathcal{P}}(\mathcal{H}_0,\mathcal{H})$ as the
 smallest integer such that every edge coloring of $\mathcal{H}_0$
 with at least $\phi_{\mathcal{P}}(\mathcal{H}_0,\mathcal{H})$ colors contains a copy
 of $\mathcal{H}$ that satisfies property $\mathcal{P}$.

\bpm
Study the problems analogous to those investigated above,
 on the more general class of uniform hypergraphs.
\epm

\bsk

\paragraph{Acknowledgements.}

This research was supported in part by the
 National Research, Development and Innovation Office,
 NKFIH Grant FK 132060.

\bsk


\begin{thebibliography}{99}
\label{s:biblio}
\addcontentsline{toc}{section}{Bibliography}

\bibitem{AC-24}
M. Axenovich, F. C. Clemen:
Rainbow subgraphs in edge-colored complete graphs: 
   Answering two questions by Erd\H os and Tuza.
J. Graph Theory 106:1 (2024), 57--66.

\bibitem{BGR-15}
A. Bialostocki, S. Gilboa, Y. Roditty:
Anti-Ramsey numbers of small graphs.
Ars Combinatoria 123 (2015), 41--53.

\bibitem{CPS-22}
Y. Caro, M. Petru\v sevski, R. \v Skrekovski:
Remarks on odd colorings of graphs.
Discrete Applied Math. 321 (2022), 392--401.

\bibitem{CPS-23}
Y. Caro, M. Petru\v sevski, R. \v Skrekovski:
Remarks on proper conflict-free colorings of graphs.
Discrete Math. 346 (2023), paper 113221.

\bibitem{ar-G}
Y. Caro, Zs. Tuza:
Monochromatic graph decompositions inspired by anti-Ramsey colorings.
arXiv:2405.19812

\bibitem{ar-nonham}
Y. Caro, Zs. Tuza:
Rainbow coloring: General theory and growth order.
Manuscript in preparation, 2024.

\bibitem{C+al-23+}
L. S. Chandran, T. Hashim, D. Jacob, R. Mathew, D. Rajendraprasad, N. Singh:
New bounds on the anti-Ramsey numbers of star graphs.
arXiv:1810.00624v2

\bibitem{ER-50}
P. Erd\H os, R. Rado:
A combinatorial theorem.
J. London Math. Soc. 25 (1950), 249--255.

\bibitem{ES-66}
P. Erd\H os, M. Simonovits:
 A limit theorem in graph theory.
Studia Sci. Math. Hungar. 1 (1966), 51--57.

\bibitem{ESS-75}
P. Erd\H os, M. Simonovits, V. T. S\'os:
 Anti-Ramsey theorems.
In: Infinite and Finite Sets,
Proc. colloq. dedicated to P. Erd\H os on his 60th birthday,
 Keszthely, 1973;
  Colloq. Math. Soc. J\'anos Bolyai, 10, Vol. II, 633--643, North-Holland, Amsterdam, 1975.

\bibitem{ES-46}
P. Erd\H os, A. H. Stone:
 On the structure of linear graphs,
Bull. Amer. Math. Soc. 52 (1946), 1087--1091.

\bibitem{dyn-surv}
S. Fujita, C. Magnant, Y. Mao, K. Ozeki:
Rainbow generalizations of Ramsey theory -- A dynamic survey.
Theory and Applications of Graphs Vol. 0: Iss. 1 (2014)
DOI: 10.20429/tag.2014.000101 
https://digitalcommons.georgiasouthern.edu/tag/vol0/iss1/1

\bibitem{GH-17}
S. Gilboa, D. Hefetz:
On degree anti-Ramsey numbers.
European J. Combinatorics 60 (2017), 31--41.

\bibitem{GR-16}
S. Gilboa, Y. Roditty:
Anti-Ramsey numbers of graphs with small connected components.
Graphs Combin. 32 (2016), 649--662.

\bibitem{GH-22+}
I. Goorevitch, R. Holzman:
Rainbow triangles in families of triangles.
arXiv:2209.15493v2

\bibitem{H-24t}
I. Harris:
Avoiding $k$-Rainbow Graphs in Edge Colorings of $K_n$
 and other Families of Graphs.
PhD Thesis, Auburn university, AL, 2024.
https://etd.auburn.edu/bitstream/handle/10415/9200/I\%20Harris\%20
Dissertation\%20Final.pdf?sequence=2

(accessed on May 8, 2024)

\bibitem{J-02}
T. Jiang:
Edge-colorings with no large polychromatic stars.
Graphs Combin. 18:2 (2002), 303--308.

\bibitem{MSTV-96}
Y. Manoussakis, M. Spyratos, Zs. Tuza, M. Voigt:
Minimal colorings for properly colored subgraphs.
Graphs Combin. 12:4 (1996), 345--360.

\bibitem{MB-06}
J. J. Montellano-Ballesteros:
On totally multicolored stars.
J. Graph Theory 51:3 (2006), 225--243.

\bibitem{MBNL-05}
J. J. Montellano-Ballesteros, V. Neumann-Lara:
An anti-Ramsey theorem on cycles. 
Graphs Combin. 21:3 (2005), 343--354.

\bibitem{P-91}
L. Pyber:
Covering the edges of a graph by...
In: Sets, Graphs and Numbers, Colloquia
Mathematica Societatis J\'anos Bolyai 60 (1991), 583--610.

\bibitem{S-04}
I. Schiermeyer:
Rainbow numbers for matchings and complete graphs.
Discrete Mathematics 286 (2004), 157--162.

\bibitem{WBZL-24}
F. Wu, H. Broersma, S. Zhang, B. Li:
Properly colored and rainbow $C_4$'s in edge-colored graphs.
J. Graph Theory 105:1 (2024), 110--135.

\bibitem{XLL-21}
J. Xu, M. Lu, K. Li:
Anti-Ramsey problems for cycles.
Applied Mathematics and Computation 408 (2021), paper 126345.

\bibitem{Y-21+}
L. T. Yuan:
The anti-Ramsey number for paths.
arXiv: 2102.00807

\end{thebibliography}
\end{document}